\documentclass[a4paper,10pt]{article}
\usepackage[utf8]{inputenc}
\usepackage[T1]{fontenc}
\usepackage{amsmath}
\usepackage{graphicx}
\usepackage{float}
\usepackage{amsfonts,amssymb}
\usepackage[numbers]{natbib}
\usepackage[export]{adjustbox}
\usepackage{subcaption}
\usepackage[justification=centering]{caption}
\DeclareMathOperator*{\argmin}{arg\,min}

\date{}
\title{\textbf{Macrophages trajectories smoothing by evolving curves}}
\author{Giulia Lupi, Karol Mikula, Seol Ah Park}
\begin{document}
\maketitle
\abstract{When analyzing cell trajectories, we often have to deal with noisy data due to the random motion of the cells and possible imperfections in cell center detection. To smooth these trajectories, we present a mathematical model and numerical method based on evolving open-plane curve approach in the Lagrangian formulation. The model contains two terms: the first is the smoothing term given by the influence of local curvature, while the other attracts the curve to the original trajectory. We use the flowing finite volume method to discretize the advection-diffusion partial differential equation. The PDE includes the asymptotically uniform tangential redistribution of curve grid points. We present results for macrophage trajectory smoothing and define a method to compute the  cell velocity for the discrete points on the smoothed curve.}
\section{Mathematical model}
In this paper, we present a mathematical model and numerical method for macrophage trajectories smoothing based on evolving open plane curve approach in the Lagrangian formulation, namely, we solve the equation
\begin{equation}
 \frac{\partial \textbf{x}}{\partial t}=-\delta k\textbf{N}+\lambda[(\textbf{x}_0-\textbf{x})\cdot \textbf{N}]\textbf{N}+\alpha \textbf{T}, 
 \label{prima}
\end{equation}
where \textbf{x} is any point of the evolving curve. The terms evolving the curve in the normal direction \textbf{N} and the
tangential direction \textbf{T} are given in details below.\\
The paper is structured as follows: in this section, we discuss the mathematical model and its possible applications. In section \ref{S2}, we show the numerical discretization of the suggested model. In section \ref{S3}, we expose the numerical experiments made on cell trajectories. Moreover, the smoothed curves allow us to compute the velocity of the cell. To do that, we propose a method to find the new velocity for the discrete points on the evolving curve based on the evolution of the length of the segments composing the piecewise linear initial condition.\\\\
Let $\Gamma$ be an open plane curve with fixed endpoints
\begin{equation}
\begin{split}
 \Gamma:[0&,1]\rightarrow\mathbb{R}^2,\\
 &u\mapsto \textbf{x}(u),
 \end{split}
\end{equation}
where $\textbf{x}(u)=(x_1(u),x_2(u))$ is a point of the curve $\Gamma$.
We consider the time evolution of the curve $\Gamma$, $\Gamma_t=\{\textbf{x}(t,u), u\in [0,1]\}$ where $t$ is the time and $\textbf{x}(t,u)$ represents the position of the point $\textbf{x}(u)\in\Gamma_t$ at time $t$. From now on, we will write $\textbf{x}$ instead of $\textbf{x}(t,u)$ for clarity.\\ 
The evolution of a point $\textbf{x}\in\Gamma$ is driven by the following general equation
\begin{equation}
 \frac{\partial \textbf{x}}{\partial t}=\textbf{V}(\textbf{x},t),
 \label{EQ1.1}
\end{equation}
where $\frac{\partial \textbf{x}}{\partial t}$ is the time derivative representing the velocity of the evolving curve. The velocity vector field $\textbf{V}$ that drives the curve motion is given by
\begin{equation}
 \textbf{V}(\textbf{x},t)=-\delta k(\textbf{x},t)\textbf{N}(\textbf{x},t)+\lambda[(\textbf{x}_0-\textbf{x})\cdot \textbf{N}(\textbf{x},t)]\textbf{N}(\textbf{x},t),
 \label{EQ1.2}
\end{equation}
where $\textbf{N}$ is the unit normal vector to the evolving curve, $k$ is its curvature, $\textbf{x}_0$ is the initial condition, " $\cdot$ " is the scalar product and $\delta,~\lambda$ are given positive parameters. From now on, we will indicate by
\begin{equation}
 w(\textbf{x},t)=(\textbf{x}_0-\textbf{x})\cdot \textbf{N}(\textbf{x},t).
 \label{w}
\end{equation}
The curvature term $-k\textbf{N}$, weighted by the parameter $\delta$, regularizes the curve while $w$, weighted by the parameter $\lambda$, attracts the evolving curve next to the initial condition. \\
The function $(\textbf{x}_0-\textbf{x})$ is the vector that realizes the minimum distance between the point $\textbf{x}$ and the initial curve $\textbf{x}_0$. We define it in the following way: 
fix a point $\textbf{x}=\textbf{x}(u)\in \Gamma_t$ on the evolving curve $\Gamma_t$, parametrized by $u\in [0,1]$. Then define
\begin{equation}
\begin{split}
 (\textbf{x}_0-\textbf{x})(u)=\argmin_{\textbf{v}\in \chi_u}\lvert \textbf{v}\rvert,
 \end{split}
 \label{EQ1.3}
\end{equation}
where $\chi_u=\{\textbf{v}\in\mathbb{R}^2: \textbf{v}=\textbf{x}_0(q)-\textbf{x}(u),q \in [0,1]\}$. For every point on the evolving curve, parametrized by $u\in [0,1]$, $(\textbf{x}_0-\textbf{x})$ selects the point on the original curve, parametrized by $q\in[0,1]$, which realizes the shortest distance and considers the vector $\textbf{v}=\textbf{x}_0(q)-\textbf{x}(u)$. Then, since the tangential velocity does not influence the shape of the evolving curve, we consider only $w(\textbf{x},t)$, the component of the vector in the normal direction. \\
We rewrite equation (\ref{EQ1.1}) as
\begin{equation}
 \frac{\partial \textbf{x}}{\partial t}=\beta\textbf{N},
 \label{EQ1.4}
\end{equation}
where the normal velocity $\beta$ is given by
\begin{equation}
 \beta=-\delta k+\lambda w.
 \label{EQ1.5}
\end{equation}
On the other hand, from the numerical point of view, the tangential velocity redistributes the points on the curve (see \cite{hou1994removing,kimura1997numerical,mikula2004direct,mikula2008simple}): this improves both the stability and the robustness of numerical computations. To redistribute the points asymptotically uniformly, we enrich the model as follows
\begin{equation}
  \frac{\partial \textbf{x}}{\partial t}=\beta\textbf{N}+\alpha \textbf{T},
  \label{EQ1.6}
\end{equation}
where $\beta$ is defined in (\ref{EQ1.5}) and $\alpha$ will be defined in section \ref{S2.1}.\\\\
Equation (\ref{prima}) can be helpful in different applications where one needs to smooth a curve and remove the noise while keeping the curve similar to the original. 
 In \cite{henry2013phagosight,park2020macrophage,Sora}, models have been proposed for segmentation and tracking of immune system cells. Once these trajectories are obtained, one may be interested in finding the actual velocity of the cells (see for example \cite{kadirkamanathan2012neutrophil}). We applied the proposed model for smoothing macrophage trajectories. Our goal was to derive the actual velocity of the cell starting from the curve found connecting the cell center in every frame of a $2$D video \cite{Sora}. Consequently, we needed to remove the random motion of the cell and correct possible imperfections in cell centers detection while keeping the final trajectory as close to the original as possible. We discuss the results in Section \ref{S3}.

\section{Numerical Algorithm}
\label{S2}
To apply the numerical scheme, the curve is discretized to a set of points, $\textbf{x}_0,~\textbf{x}_1,...,~\textbf{x}_{n+1}$ as displayed in Fig. \ref{Fig1}.\\
Let $\lvert\textbf{x}_u\rvert>0$, where $\textbf{x}_u=(\frac{\partial x_1}{\partial u},\frac{\partial x_2}{\partial u})$ and $g=\lvert\textbf{x}_u\rvert=\sqrt{(\frac{\partial x_1}{\partial u})^2+(\frac{\partial x_2}{\partial u})^2}$. If we denote by $s$ the unit arc-length parametrization of the curve $\Gamma$, then $ds=\lvert\textbf{x}_u\rvert du=g du$ and $du=\frac{1}{g}ds$. The unit tangent vector $\textbf{T}$ is defined as $\textbf{T}=\frac{\partial \textbf{x}}{\partial s}=\textbf{x}_s$ while the unit normal vector $\textbf{N}$ is $\textbf{N}=\textbf{x}_{s}^{\perp}$ such that $\textbf{T}\wedge\textbf{N}=-1$. If $\textbf{T}=(\frac{x_1}{\partial s},\frac{x_2}{\partial s})$ then $\textbf{N}=(\frac{x_2}{\partial s},-\frac{x_1}{\partial s})$. From the Frenet-Serret formulas we have $\textbf{T}_s=-k\textbf{N}$ and $\textbf{N}_s=k\textbf{T}$, where $k$ is the curvature. It follows that $-k\textbf{N}=\textbf{T}_s=(\textbf{x}_s)_s=\textbf{x}_{ss}$.\\
We can rewrite (\ref{EQ1.6}) into the form of the so-called intrinsic partial differential equation
\begin{equation}
 \textbf{x}_t=\delta\textbf{x}_{ss}+\alpha\textbf{x}_s+\lambda w\textbf{x}_{s}^{\perp}
 \label{EQ1.7}
\end{equation}
which is suitable for numerical discretization. Since $\textbf{x}=(x_1,x_2)$, (\ref{EQ1.7}) represents a system of two partial differential equations for the two components of the position vector $\textbf{x}$. 
\subsection{Suitable choice of tangential velocity}
In \cite{hou1994removing,kimura1997numerical,mikula2004direct,mikula2008simple},  it has been shown how to choose the tangential velocity to have an asymptotically uniform redistribution of points; in this section, we recall it for completeness. \\
If we want to redistribute the points along the curve, we have to consider the ratio $\frac{g}{L}$, where $g=\lvert \partial_u\textbf{x}\rvert$ represents the local length of the curve while $L$ is the global length. In the discrete form, we obtain
\[
 \frac{g}{L}\approx \frac{\frac{\lvert \textbf{x}_i-\textbf{x}_{i-1}\rvert}{h}}{L}=\frac{\lvert \textbf{x}_i-\textbf{x}_{i-1}\rvert}{Lh}=\frac{\lvert \textbf{x}_i-\textbf{x}_{i-1}\rvert}{\frac{L}{n+1}},
\]
where $h=\frac{1}{n+1}$. We define the $i$-th element $\textbf{h}_i=\textbf{x}_i-\textbf{x}_{i-1}$ with length $h_i=\lvert \textbf{h}_i\rvert$. In the formula above, $n+1$ is the number of elements used in the spatial discretization based on the flowing finite volume method. The numerator represents the distance between two grid points; on the other hand, the denominator represents the distance that would occur if we had uniformly distributed points. It is then clear that if we want uniformly distributed points, we will have to ask for that ratio to tend to $1$. In continuous setting, we have to fulfill the condition $\frac{g}{L}\rightarrow 1$.\\
For the time evolution of the ratio, we obtain
\begin{equation}
 (\frac{g}{L})_t=\frac{g_t L-L_t g}{L^2}.
 \label{EQ1.8}
\end{equation}
We observe first that the time derivative of the local length is obtained using the Frenet-Serret formulas and is given by
\begin{equation}
 \begin{split}
  g_t&=\lvert\textbf{x}_u\rvert_t=\frac{\textbf{x}_u}{\lvert\textbf{x}_u\rvert}\cdot(\textbf{x}_u)_t=\frac{g\textbf{x}_s}{g}\cdot(\textbf{x}_t)_u=\textbf{T}\cdot(\beta \textbf{N}+\alpha\textbf{T})_u=\textbf{T}\cdot g(\beta\textbf{N}+\alpha\textbf{T})_s=\\
  &=\textbf{T}\cdot g(\beta_s\textbf{N}+\beta\textbf{N}_s+\alpha_s\textbf{T}+\alpha\textbf{T}_s)=\textbf{T}\cdot g(\beta_s\textbf{N}+k\beta\textbf{T}+\alpha_s\textbf{T}-k\alpha\textbf{N})=\\
  &=gk\beta+g\alpha_s=gk\beta+\alpha_u.
 \end{split}
 \label{EQ1.9}
\end{equation}
To get the formula for the global length, we have to integrate $g=\lvert\textbf{x}_u\rvert$ over the curve $\Gamma$
\[
 L=\int_{\Gamma}g du=\int_{0}^{1} g du,
\]
and then
\begin{equation}
 L_t=\int_{0}^{1}g_tdu=\int_{0}^{1}gk\beta du+\int_{0}^{1} \alpha_u du=\int_{0}^{1} k\beta ds+\alpha(1)-\alpha(0).
 \label{EQ1.10}
\end{equation}
Now, since we fixed the endpoints of the open curve, we get $\alpha(1)=\alpha(0)=0$, and equation (\ref{EQ1.10}) reduces to 
\begin{equation}
 L_t=\int_{0}^{1} k\beta ds=L\langle k\beta\rangle_{\Gamma},
 \label{EQ1.11}
\end{equation}
where $\langle k\beta\rangle_{\Gamma}=\frac{1}{L}\int_{\Gamma}k\beta ds$.\\
Substituting (\ref{EQ1.9}) and (\ref{EQ1.11}) in (\ref{EQ1.8}) we obtain
\begin{equation}
 (\frac{g}{L})_t=\frac{(gk\beta+\alpha_u)L-Lg\langle k\beta\rangle_{\Gamma}}{L^2}=\frac{g}{L}(k\beta+\alpha_s-\langle k\beta\rangle_{\Gamma}).
\end{equation}
If we now impose that $(\frac{g}{L})_t=\omega(1-\frac{g}{L})$, where $\omega$ is a parameter that determines the speed of convergence of the ratio, we get the desired condition $\frac{g}{L}\rightarrow 1$. Finally, we obtain the formula for the tangential velocity 
\begin{equation}
 \alpha_s=\langle k\beta\rangle_{\Gamma}-k\beta+\omega(\frac{L}{g}-1).
 \label{EQ2.1.14}
\end{equation}

\label{S2.1}
\subsection{Numerical Discretization}
\label{S2.2}
Let us consider the following form of the intrinsic PDE (\ref{EQ1.7})
\begin{equation}
 \textbf{x}_t-\alpha\textbf{x}_s=\delta\textbf{x}_{ss}+\lambda w\textbf{x}_{s}^{\perp},
 \label{EQ1.2.1}
\end{equation}
where $\alpha$ and $w$ are given by (\ref{EQ2.1.14}) and (\ref{w}) respectively. First, we perform the spatial discretization based on the flowing finite volume method \cite{sevcovic2001evolution,mikula2008simple}, then we discuss the semi-implicit time discretization, which is implicit in the intrinsic diffusion term and uses the inflow-implicit/outflow-explicit strategy for intrinsic advection term \cite{mikula2011inflow}.\\
 \begin{figure} 
  \includegraphics[width=\linewidth]{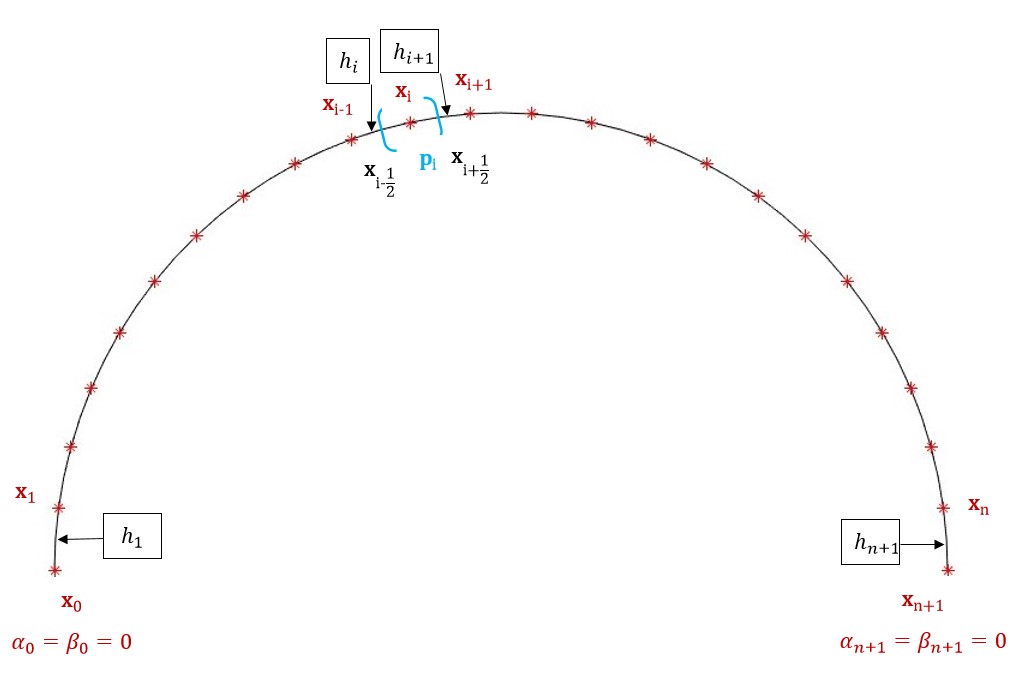} \\
  \caption{Open curve discretization. Finite volume $\textbf{p}_i=[\textbf{x}_{i-\frac{1}{2}}, \textbf{x}_{i+\frac{1}{2}}]$ is highlighted in blue.}
  \label{Fig1}
\end{figure} 
Consider an open curve discretized into $n+2$ grid points. We fixed the two points $ \textbf{x}_0$ and $ \textbf{x}_{n+1}$ at the endpoints of the curve. Consequently, for these two points, the velocity is zero.\\ 
We define $h_i=\lvert  \textbf{x}_i- \textbf{x}_{i-1}\rvert$ and consider the finite volume $\textbf{p}_i=[ \textbf{x}_{i-\frac{1}{2}}, \textbf{x}_{i+\frac{1}{2}}]$ where $ \textbf{x}_{i+\frac{1}{2}}$ represents the middle point between $ \textbf{x}_i$ and $ \textbf{x}_{i+1}$, i.e.,
\[
  \textbf{x}_{i+\frac{1}{2}}=\frac{ \textbf{x}_i+ \textbf{x}_{i+1}}{2}.
\]
Notice that the velocities $\alpha_i,~\beta_i$ will be referred to the points $ \textbf{x}_i, i=0,...,n+1$ with $\alpha_0=0$, $\beta_0=0$ and $\alpha_{n+1}=0$, $\beta_{n+1}=0$. On the other side, the curvature $k_i$ will be referred to the elements $\textbf{h}_i,~i=1,...,n+1$. \\
Integrating (\ref{EQ1.2.1}) over the finite volume $\textbf{p}_i$, see Fig. \ref{Fig1}, we get
\begin{equation}
 \int_{ \textbf{x}_{i-\frac{1}{2}}}^{ \textbf{x}_{i+\frac{1}{2}}} \textbf{x}_tds -\int_{ \textbf{x}_{i-\frac{1}{2}}}^{ \textbf{x}_{i+\frac{1}{2}}} \alpha\textbf{x}_sds=\delta \int_{ \textbf{x}_{i-\frac{1}{2}}}^{ \textbf{x}_{i+\frac{1}{2}}} \textbf{x}_{ss}ds+\lambda\int_{ \textbf{x}_{i-\frac{1}{2}}}^{ \textbf{x}_{i+\frac{1}{2}}}w\textbf{x}_{s}^{\perp}ds,
 \label{EQ1.2.3n}
\end{equation}
where $\delta,~\lambda$ are constants and the values $\alpha,w$ are considered to be constant over the finite volume $\textbf{p}_i$ and will be indicated as $\alpha_i,~w_i$. \\
The length of the interval $\textbf{p}_i$ is $\frac{h_i+h_{i+1}}{2}$ consequently, using the Newton-Leibniz formula we get the following approximation of (\ref{EQ1.2.3n}) for $i=1,...,n$
\begin{equation}
 \frac{h_i+h_{i+1}}{2}(\textbf{x}_i)_t-\alpha_i[\textbf{x}]_{ \textbf{x}_{i-\frac{1}{2}}}^{ \textbf{x}_{i+\frac{1}{2}}}=\delta[\textbf{x}_s]_{ \textbf{x}_{i-\frac{1}{2}}}^{ \textbf{x}_{i+\frac{1}{2}}}+\lambda w_i([\textbf{x}]_{ \textbf{x}_{i-\frac{1}{2}}}^{ \textbf{x}_{i+\frac{1}{2}}})^{\perp}.
 \label{EQ1.2.3}
\end{equation}
First of all, notice that
\begin{align}
 -\alpha_i[\textbf{x}]_{ \textbf{x}_{i-\frac{1}{2}}}^{ \textbf{x}_{i+\frac{1}{2}}}&=-\alpha_i( \textbf{x}_{i+\frac{1}{2}}- \textbf{x}_{i-\frac{1}{2}})=-\alpha_i(\frac{\textbf{x}_i+\textbf{x}_{i+1}}{2}-\frac{\textbf{x}_i+\textbf{x}_{i-1}}{2})=\\
 &=-\alpha_i(\frac{\textbf{x}_{i+1}-\textbf{x}_{i-1}}{2})=\frac{\alpha_i}{2}(\textbf{x}_i-\textbf{x}_{i+1})-\frac{\alpha_i}{2}(\textbf{x}_i-\textbf{x}_{i-1}).
 \label{EQ1.2.4}
\end{align}
Then, approximating $\textbf{x}_s$ by a finite difference, we get
\begin{equation}
[\textbf{x}_s]_{ \textbf{x}_{i-\frac{1}{2}}}^{ \textbf{x}_{i+\frac{1}{2}}}=(\frac{\textbf{x}_{i+1}-\textbf{x}_i}{h_{i+1}}-\frac{\textbf{x}_{i}-\textbf{x}_{i-1}}{h_{i}}).
\label{EQ1.2.6}
\end{equation}
Combining (\ref{EQ1.2.4}), (\ref{EQ1.2.6}), and (\ref{EQ1.2.3}), we obtain
\begin{equation}
\begin{split}
 &\frac{h_i+h_{i+1}}{2}(\textbf{x}_i)_t+\frac{\alpha_i}{2}(\textbf{x}_i-\textbf{x}_{i+1})-\frac{\alpha_i}{2}(\textbf{x}_i-\textbf{x}_{i-1})=\\
 &=\delta_i(\frac{\textbf{x}_{i+1}-\textbf{x}_i}{h_{i+1}}-\frac{\textbf{x}_{i}-\textbf{x}_{i-1}}{h_{i}})+\lambda w_i(\frac{\textbf{x}_{i+1}-\textbf{x}_{i-1}}{2})^{\perp}.
 \label{EQ1.2.7}
 \end{split}
\end{equation}
The basic idea of the inflow-implicit/outflow-explicit strategy is that we treat outflow from a cell explicitly while inflow implicitly \cite{mikula2011inflow}. If $\alpha_i<0$, i.e., $(-\alpha)$ in (\ref{EQ1.2.1}) is positive in the finite volume $\textbf{p}_i$, then there is an inflow into the finite volume through its boundary point $\textbf{x}_{i-\frac{1}{2}}$ and there is an outflow through the other boundary point $\textbf{x}_{i+\frac{1}{2}}$. On the other hand, if $\alpha_i>0$, there is outflow through $\textbf{x}_{i-\frac{1}{2}}$ and inflow through $\textbf{x}_{i+\frac{1}{2}}$. We therefore define
\begin{equation}
\begin{split}
 b_{i-\frac{1}{2}}^{in}=max(-\alpha_i,0),b_{i-\frac{1}{2}}^{out}=min(-\alpha_i,0),\\
 b_{i+\frac{1}{2}}^{in}=max(\alpha_i,0),b_{i+\frac{1}{2}}^{out}=min(\alpha_i,0),
 \end{split}
 \label{EQ1.2.8}
\end{equation}
and rewrite equation (\ref{EQ1.2.7}) as
\begin{equation}
\begin{split}
&\frac{h_i+h_{i+1}}{2}(\textbf{x}_i)_t+\frac{1}{2}(b_{i+\frac{1}{2}}^{in}+b_{i+\frac{1}{2}}^{out})(\textbf{x}_i-\textbf{x}_{i+1})+\frac{1}{2}(b_{i-\frac{1}{2}}^{in}+b_{i-\frac{1}{2}}^{out})(\textbf{x}_i-\textbf{x}_{i-1})\\ 
&=\delta(\frac{\textbf{x}_{i+1}-\textbf{x}_i}{h_{i+1}}-\frac{\textbf{x}_{i}-\textbf{x}_{i-1}}{h_{i}})+\lambda w_i(\frac{\textbf{x}_{i+1}-\textbf{x}_{i-1}}{2})^{\perp}.
\end{split}
 \label{EQ1.2.9}
\end{equation}
Let us approximate the time derivative by the finite difference $\textbf{x}_t=\frac{\textbf{x}_{i}^{m+1}-\textbf{x}_{i}^{m}}{\tau}$: here $m$ is the time step index and $\tau$ is the length of the discrete time step. Remark that we take the unknowns in the inflow part of the advection term implicitly while we take the ones in the outflow part explicitly. Moreover, we take the diffusion term implicitly and the attracting term explicitly. \\
Finally, we obtain the fully discrete scheme
\begin{equation}
 \begin{split}
  &\textbf{x}_{i-1}^{m+1}(-\frac{\delta}{h_{i}^{m}}-\frac{b_{i-\frac{1}{2}}^{in^{m}}}{2})+\textbf{x}_{i+1}^{m+1}(-\frac{\delta}{h_{i+1}^{m}}-\frac{b_{i+\frac{1}{2}}^{in^{m}}}{2})+\\
  &\textbf{x}_{i}^{m+1}(\frac{h_{i+1}^{m}+h_{i}^{m}}{2\tau}+\frac{\delta}{h_{i}^{m}}+\frac{\delta}{h_{i+1}^{m}}+\frac{b_{i-\frac{1}{2}}^{in^{m}}}{2}+\frac{b_{i+\frac{1}{2}}^{in^{m}}}{2})=\textbf{x}_{i}^{m}\frac{h_{i+1}^{m}+h_{i}^{m}}{2\tau}\\
  &-\frac{b_{i+\frac{1}{2}}^{out^{m}}}{2}(\textbf{x}_{i}^{m}-\textbf{x}_{i+1}^{m})-\frac{b_{i-\frac{1}{2}}^{out^{m}}}{2}(\textbf{x}_{i}^{m}-\textbf{x}_{i-1}^{m})+\lambda w_{i}^{m}(\frac{\textbf{x}_{i+1}^{m}-\textbf{x}_{i-1}^{m}}{2})^{\perp},
 \end{split}
 \label{EQ1.2.10}
\end{equation}
for $i=1,...,n$ where $n$ is the number of unknown grid points.\\ 
The system (\ref{EQ1.2.10}) is represented by a strictly diagonally dominant matrix, then it is always solvable by the classical Thomas algorithm without any restriction on the time step $\tau$.\\
In the numerical scheme (\ref{EQ1.2.10}), there are two parameters $\alpha_{i}^{m}$ and $w_{i}^{m}$ given by (\ref{EQ2.1.14}) and (\ref{w}) for which we have not yet defined the discretization.\\
Let us first consider the discretization of $w_{i}^{m}$ referred to the curve grid point $\textbf{x}_{i}^{m}$. Remark that the shortest distance of a point to a curve is in the normal direction to the curve.\\
The algorithm is organized as follows: consider a point $\textbf{x}_{i}^{m}$ on the evolving curve. Set $D_{min}=D$ (where $D$ is a real number to be chosen reasonably) and for $j=1,...,n+1$ repeat: 
\begin{enumerate}
\item consider the line $r$ passing through the points $\textbf{x}_{j-1}^{0}$ and $\textbf{x}_{j}^{0}$ of the original curve,
\item find the line $s$ perpendicular to $r$ and passing through $\textbf{x}_{i}^{m}$,
\item find the point $\textbf{q}$ of intersection between the line $r$ and the line $s$, 
\item find $d=d(\textbf{q},\textbf{x}_{i}^{m})$, where $d$ is the Euclidean distance of $\textbf{q}$ and $\textbf{x}_{i}^{m}$,
\item  if $\textbf{q}\in[\textbf{x}_{j-1}^{0},\textbf{x}_{j}^{0}]$ and $d<D_{min}$ then 
\begin{equation}
\begin{split}
 &\textbf{x}^{0}-\textbf{x}^{m}_{i}=\textbf{q}-\textbf{x}_{i}^{m}\\
 &D_{min}=\lvert \textbf{x}^{0}-\textbf{x}^{m}_{i} \rvert
 \label{attracting}
 \end{split}
\end{equation}
\label{ALG}
\end{enumerate}
Due to the discretization of the curve, it may happen that the algorithm does not find any point $\textbf{q}$.
In this case, for every $j=0,...,n+1$, we check the Euclidean distance $d(\textbf{x}_{j}^{0},\textbf{x}_{i}^{m})$ and we select the minimum $j=j_{min}$.Then,
\[
 \textbf{x}^{0}-\textbf{x}^{m}_{i}=\textbf{x}_{j_{min}}^{0}-\textbf{x}_{i}^{m}.
\]

This approach is natural since we are trying to find the minimum distance between a point and a segment: it is either the minimum distance found in (\ref{attracting}) or the minimum distance between the point and one of the endpoints of the segment.\\
Therefore, we get
\begin{equation}
 w_{i}^{m}=(\textbf{x}^{0}-\textbf{x}^{m}_{i})\cdot \textbf{N}_{i}^{m},
\end{equation}
where $\textbf{N}_{i}^{m}=(\frac{\textbf{x}_{i+1}^{m}-\textbf{x}_{i-1}^{m}}{h_{i}^{m}+h_{i+1}^{m}})^{\perp}$.\\

Consider now the discretization of $\alpha_{i}^{m}$. Remark that, since we fixed the first and the last point, we set $\alpha_0=0$ and $\alpha_{n+1}=0$. We get $\alpha_{i}^{m}$ for $i=1,...,n$ by
\begin{equation}
 \alpha_{i}^{m}=\alpha_{i-1}^{m}+h_{i}^{m}\langle k\beta \rangle_{\Gamma}^{m}-h_{i}^{m}k_{i}^{m}\beta_{i}^{m}+\omega(\frac{L^m}{n+1}-h_{i}^{m}).
\end{equation}
The curvature $k_{i}^{m}$, the normal velocity $\beta_{i}^{m}$, the mean value $\langle k\beta \rangle_{\Gamma}^{m}$, and the total length $L^m$ are given by the following formulas \cite{ambroz2020semi,mikula2008simple,mikula2021automated}
\begin{equation}
 \begin{split}
  &k_{i}^{m}=sgn(\textbf{h}_{i-1}^{m}\wedge \textbf{h}_{i+1}^{m})\frac{1}{2h_{i}^{m}}arccos(\frac{\textbf{h}_{i-1}^{m}\cdot \textbf{h}_{i+1}^{m}}{h_{i-1}^{m} h_{i+1}^{m}}),\\
  &\beta_{i}^{m}=-\delta k_{i}+\lambda  w_i,\\
  &\langle k\beta \rangle_{\Gamma}^{m}=\frac{1}{L^{m}}\sum_{l=1}^{n+1}h_{l}^{m}k_{l}^{m}\beta_{i}^{m},\\
  &L^m=\sum_{l=1}^{n+1}h_{l}^{m},
  \end{split}
  \label{EQ2.30}
\end{equation}
where $\textbf{h}_{i}^{m}=\textbf{x}_{i}^{m}-\textbf{x}_{i-1}^{m}$, $\lvert\textbf{h}_{i}^{m}\rvert=h_{i}^{m}$. For the first and last elements, since we don't have the values of $h_{i-1}^m$, respectively $h_{i+1}^m$, we set $k_{1}^{m}=k_{2}^{m}$ and $k_{n+1}^{m}=k_{n}^{m}$.
\section{Numerical experiments}
\label{S3}
\subsection{Smoothed trajectories}
\label{S3.1}
\begin{figure}
\centering
  \includegraphics[width=15cm,center]{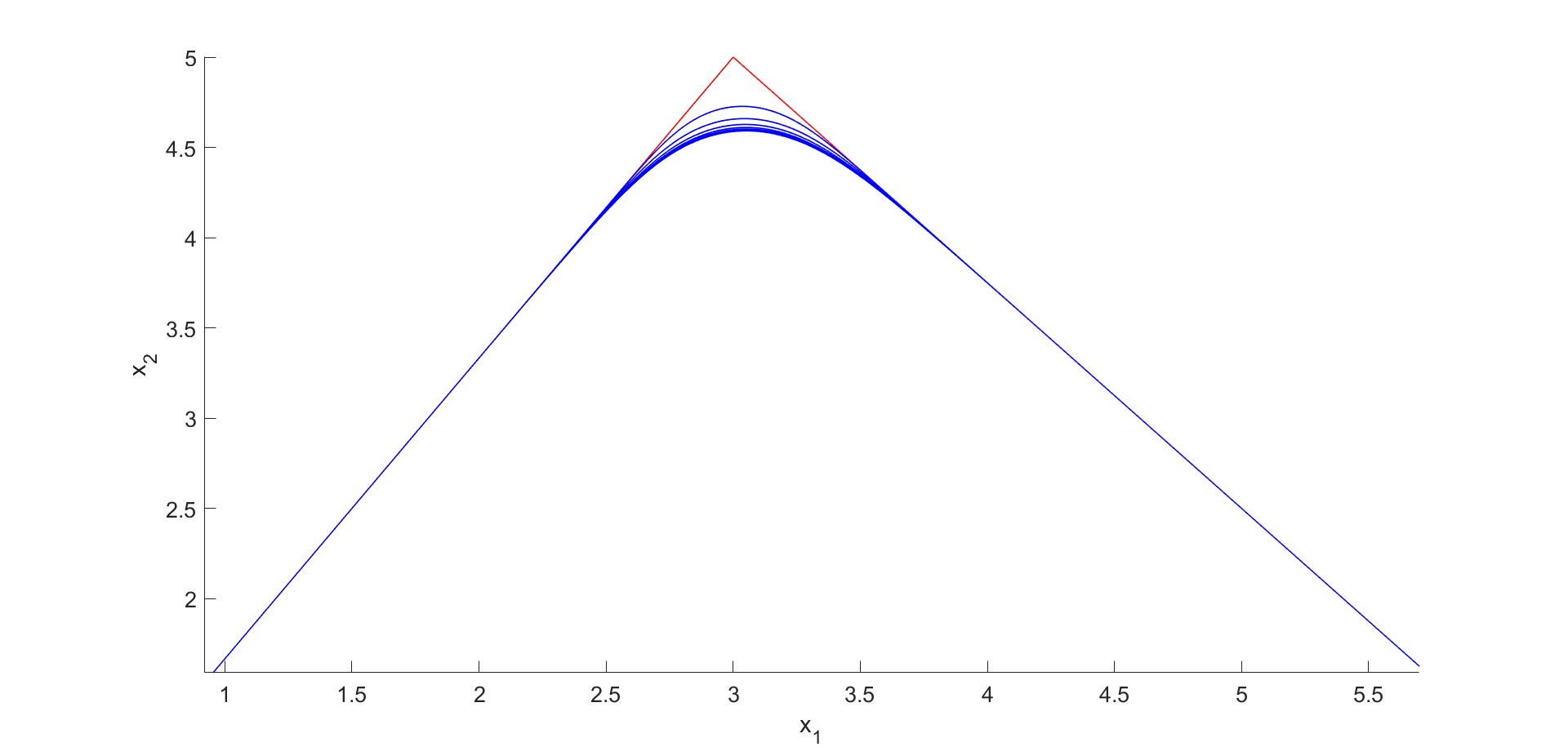}
  \includegraphics[width=15cm,center]{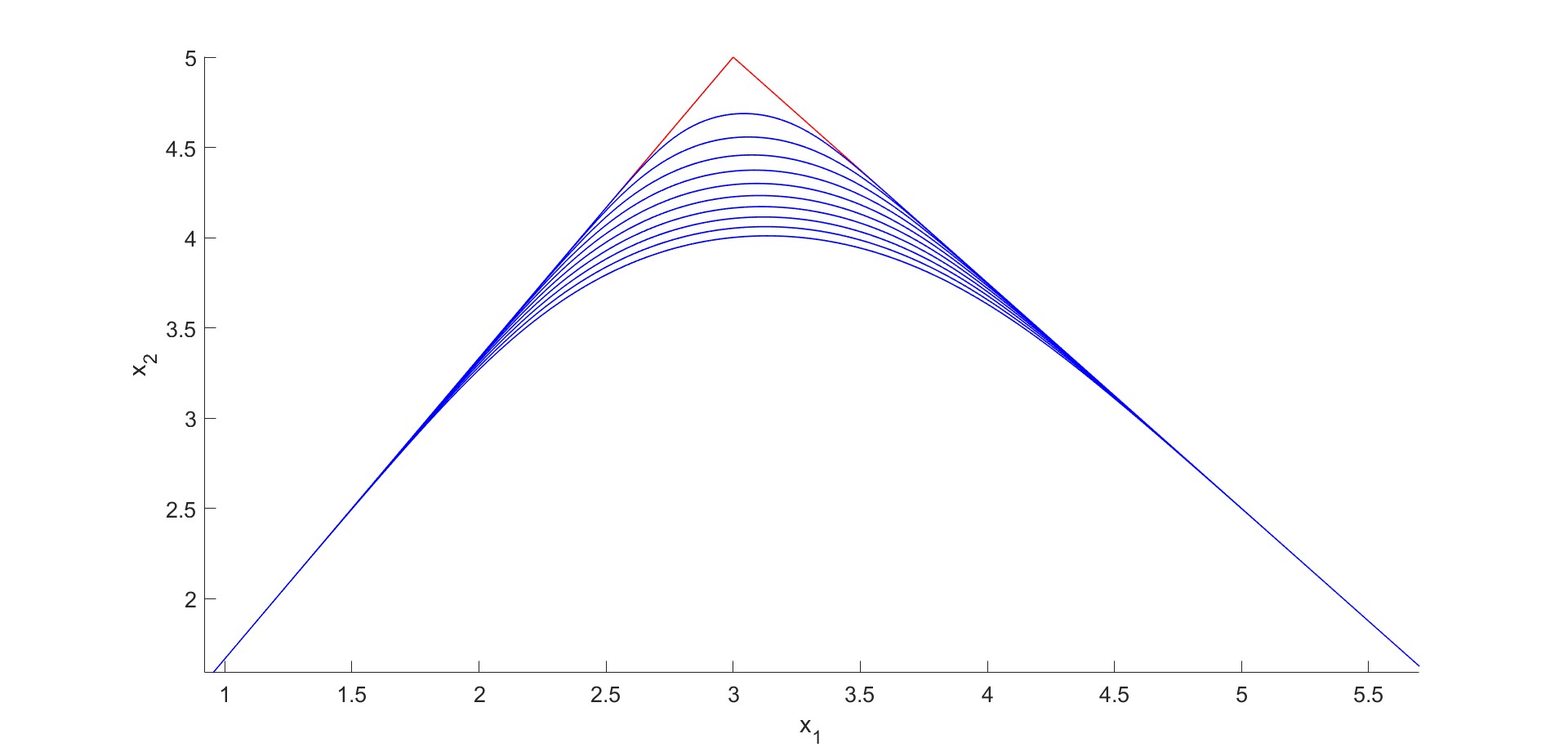}
\caption{Evolution of the initial curve (red) with $\lambda=1$  (top) and $\lambda=0$ (bottom), $\delta=0.001$, $\omega=1$, and $\tau=0.001$. The blue curves are the evolved curves plotted every $1000$ time step, from time step $1000$ to $10000$}
  \label{Fig2}
\end{figure}

\begin{figure}
\centering
  \includegraphics[width=15cm,center]{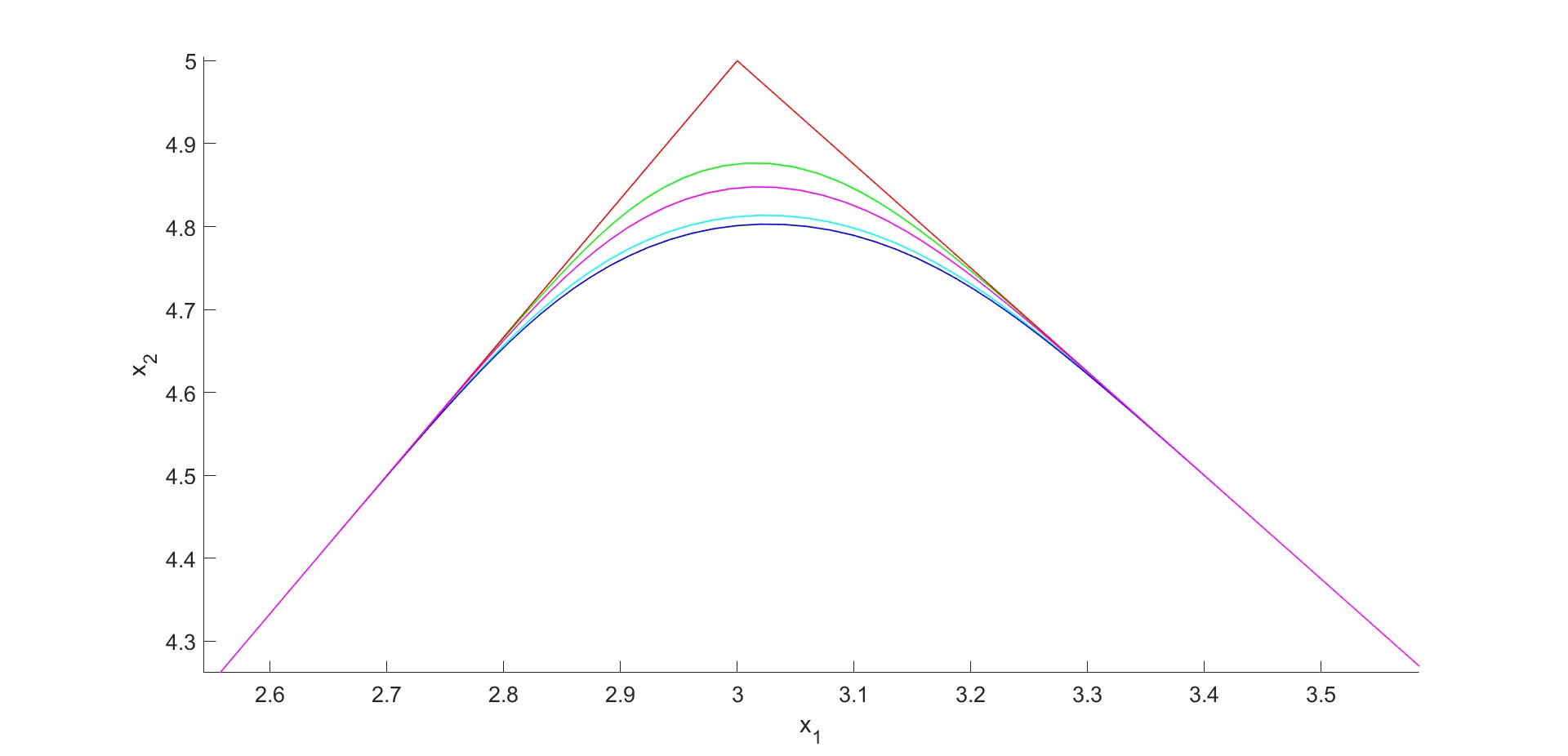}
  \caption{Comparison of the evolution of the initial curve (red) for different values of $\lambda$ after $400$ time steps. Results are shown for $\lambda=0$ (blue), $\lambda=1$ (light blue), $\lambda=5$ (pink) and $\lambda=10$ (green), $\delta=0.001$, $\omega=1$, and $\tau=0.001$.}
  \label{Figsub}
\end{figure}

\begin{figure}
  \centering
     \includegraphics[width=20cm,center]{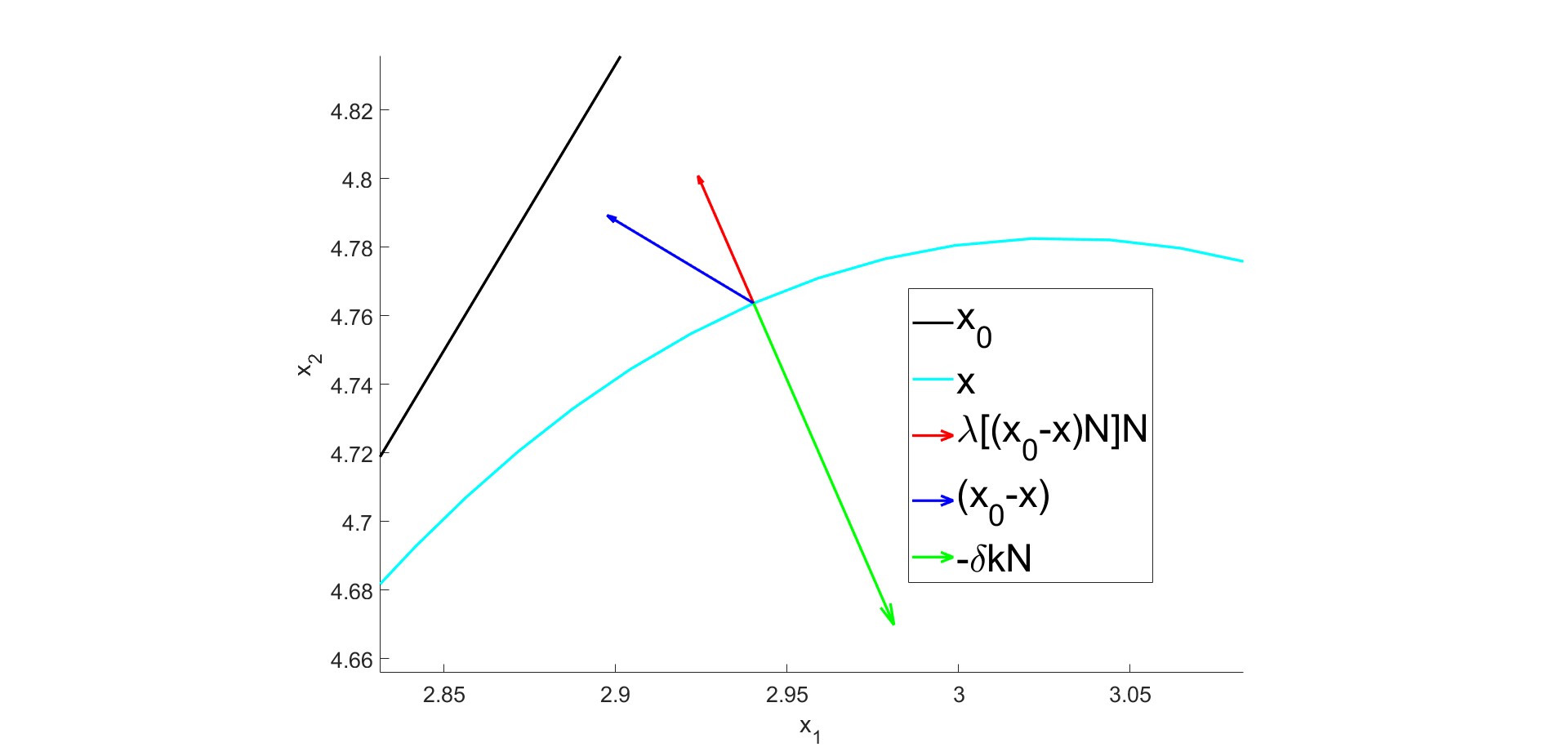}
  \caption{Plotting of vectors driving one of the curve grid points during the evolution. $-\delta k \textbf{N}$ (green), $(\textbf{x}_0-\textbf{x})$(blue),  and its projection in normal direction to the curve $[(\textbf{x}_0-\textbf{x})\cdot\textbf{N}]\textbf{N}$(red).}
\label{Figsub2}
\end{figure}

First, we considered the simple curve shown in Figs. \ref{Fig2}-\ref{Figsub2} to see the influence of the term which attracts the curve to the initial condition. In Fig. \ref{Fig2}, we considered the evolution of the curve driven by the PDE (\ref{prima}) with $\lambda=0$ and $\lambda=1$, $\delta=0.001$, $\omega=1$, and $\tau=0.001$ for $10000$ time steps. The red line represents the initial condition while the blue lines are the results of the smoothing plotted every $1000$ time step. The curve evolved only by the curvature influence goes to the straight line, while if we consider $\lambda=1$, the curve is attracted to the initial condition and stays closer to it. In Fig. \ref{Figsub} we considered the same initial curve as in Fig. \ref{Fig2} and we plotted the results of the smoothing after $400$ time steps for  $\lambda=0$, $\lambda=1$, $\lambda=5$ and $\lambda=10$, $\delta=0.001$, $\omega=1$, and $\tau=0.001$: we observe that considering bigger $\lambda$, the smoothed curve stays closer to the original one. To better understand the influence of the attracting term in Fig. \ref{Figsub2} we plotted the vectors driving one of the discrete points of the evolving curve: the blue arrow is the vector $(\textbf{x}_0-\textbf{x})$, the red arrow is its projection in the normal direction to the evolving curve, and the green one is the diffusion term. We observe that the attracting vector always points towards the original curve: it is clear that, if we consider bigger $\lambda$, the result will stay closer to the initial condition.\\

\begin{figure}
  \includegraphics[width=10cm,center]{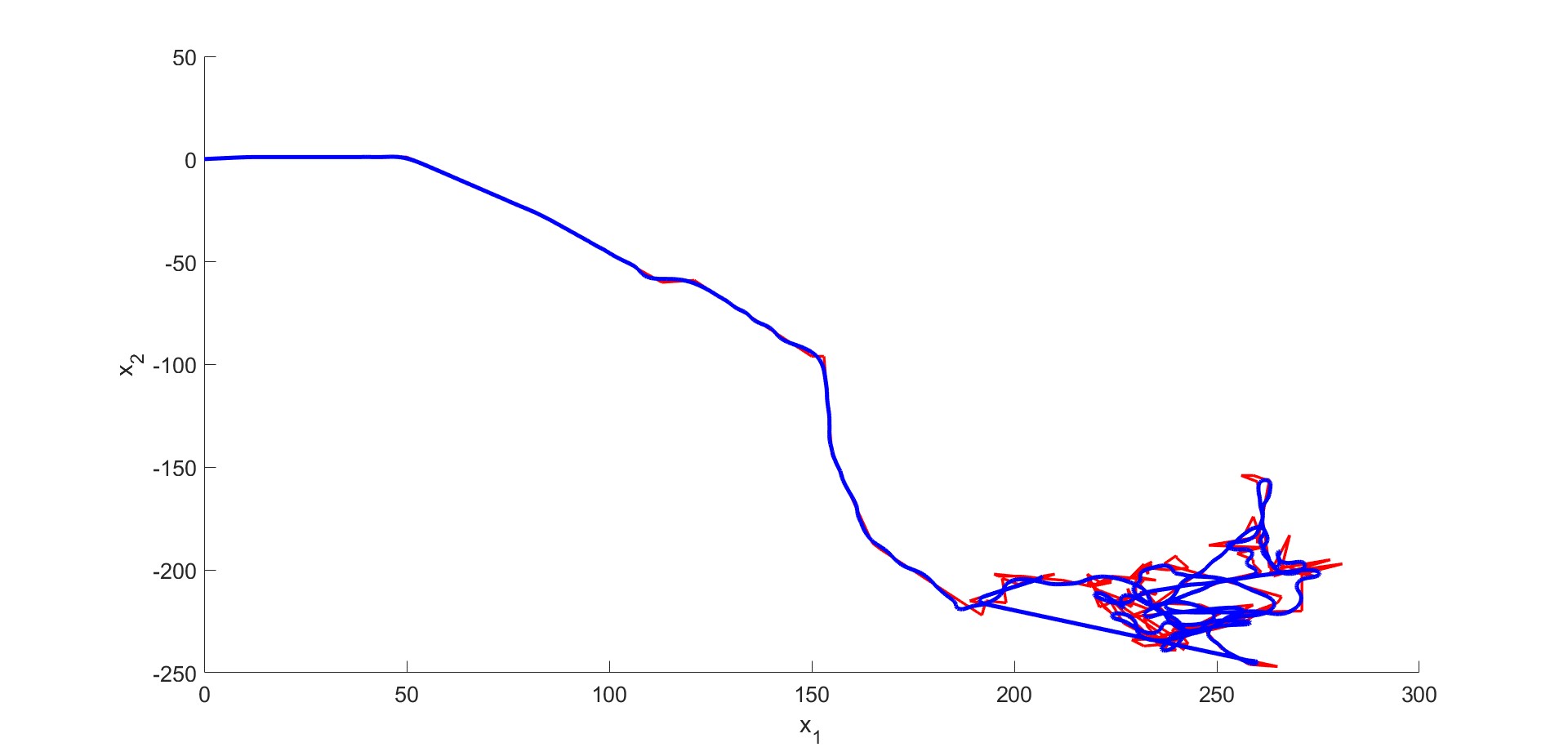}
  \includegraphics[width=10cm,center]{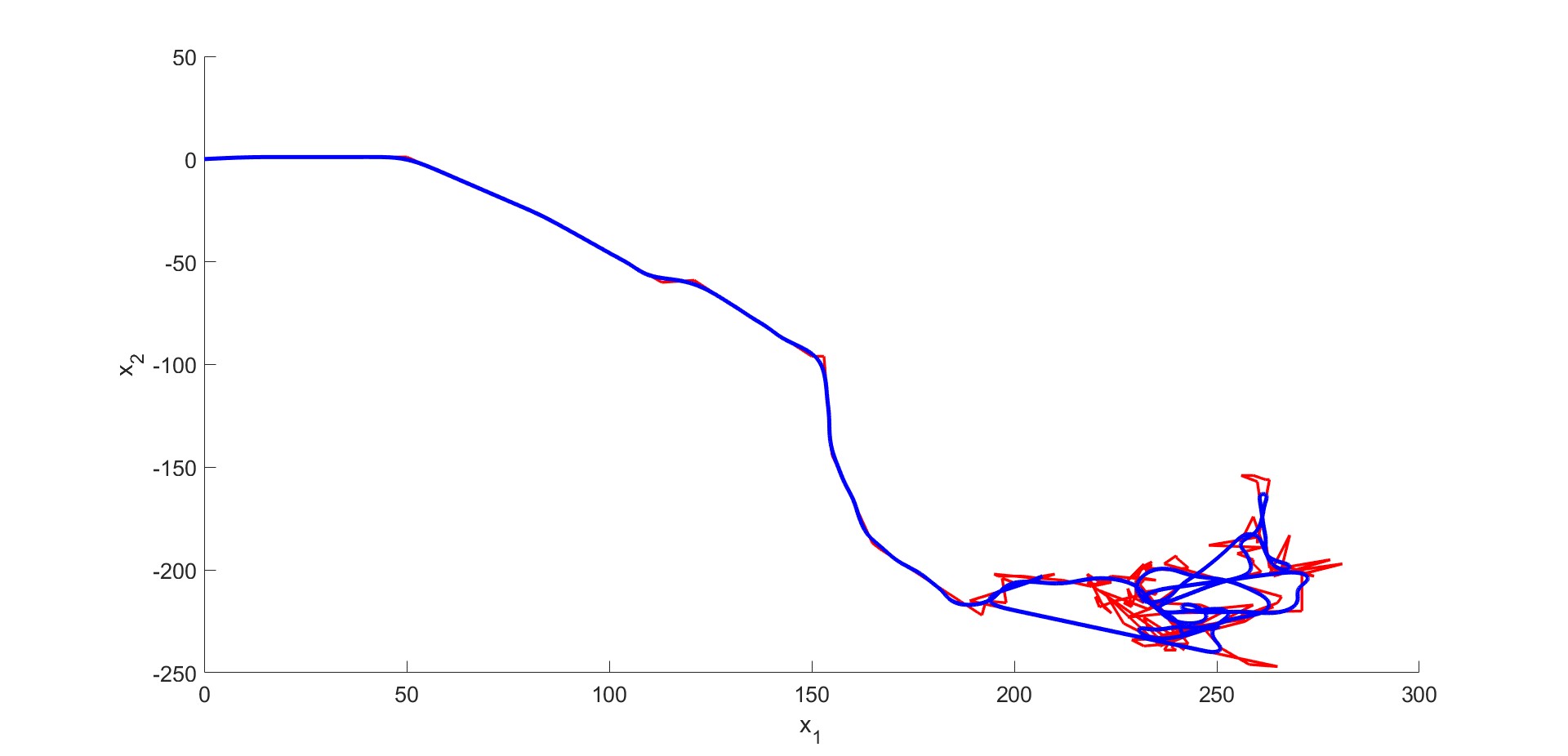}
  \includegraphics[width=10cm,center]{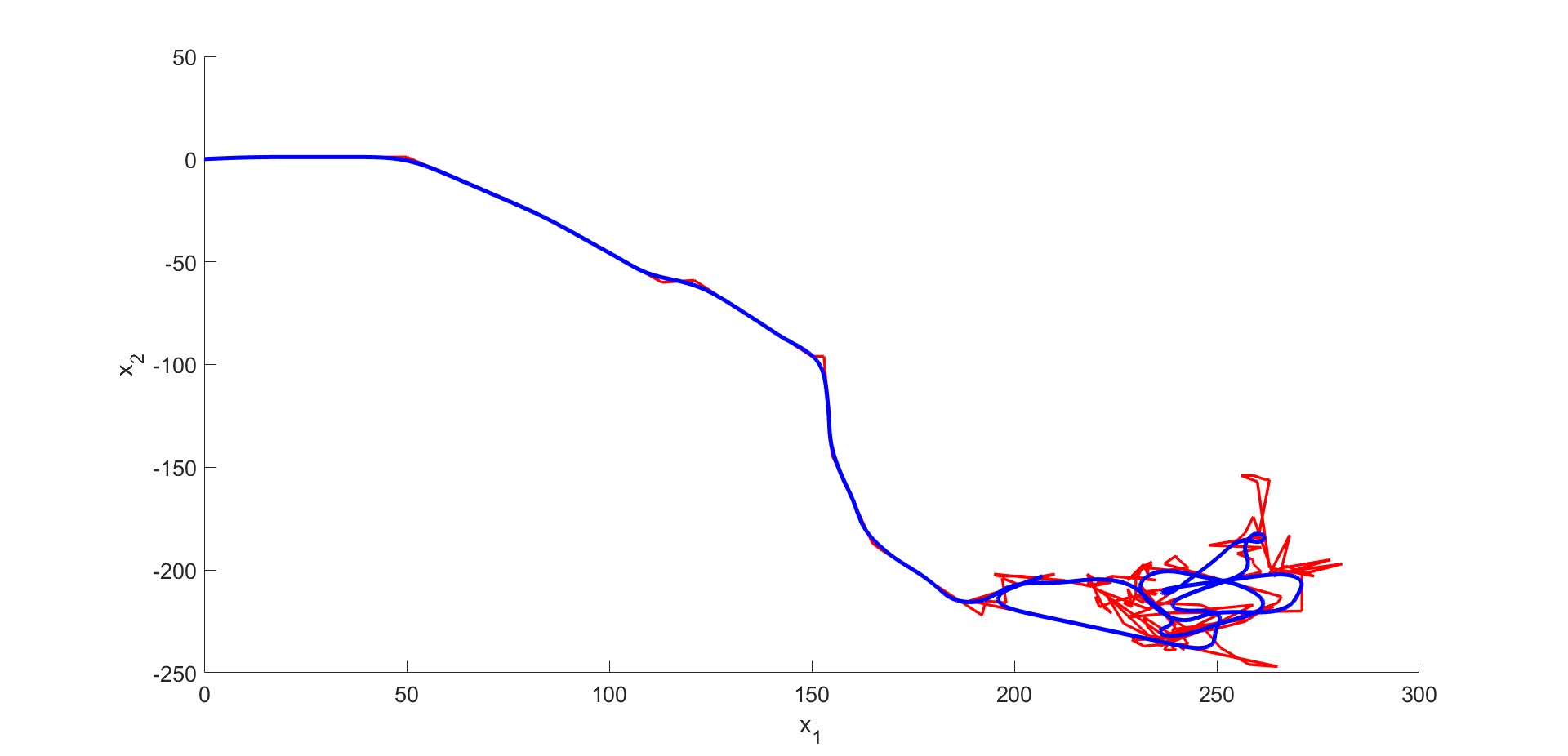}
  \includegraphics[width=10cm,center]{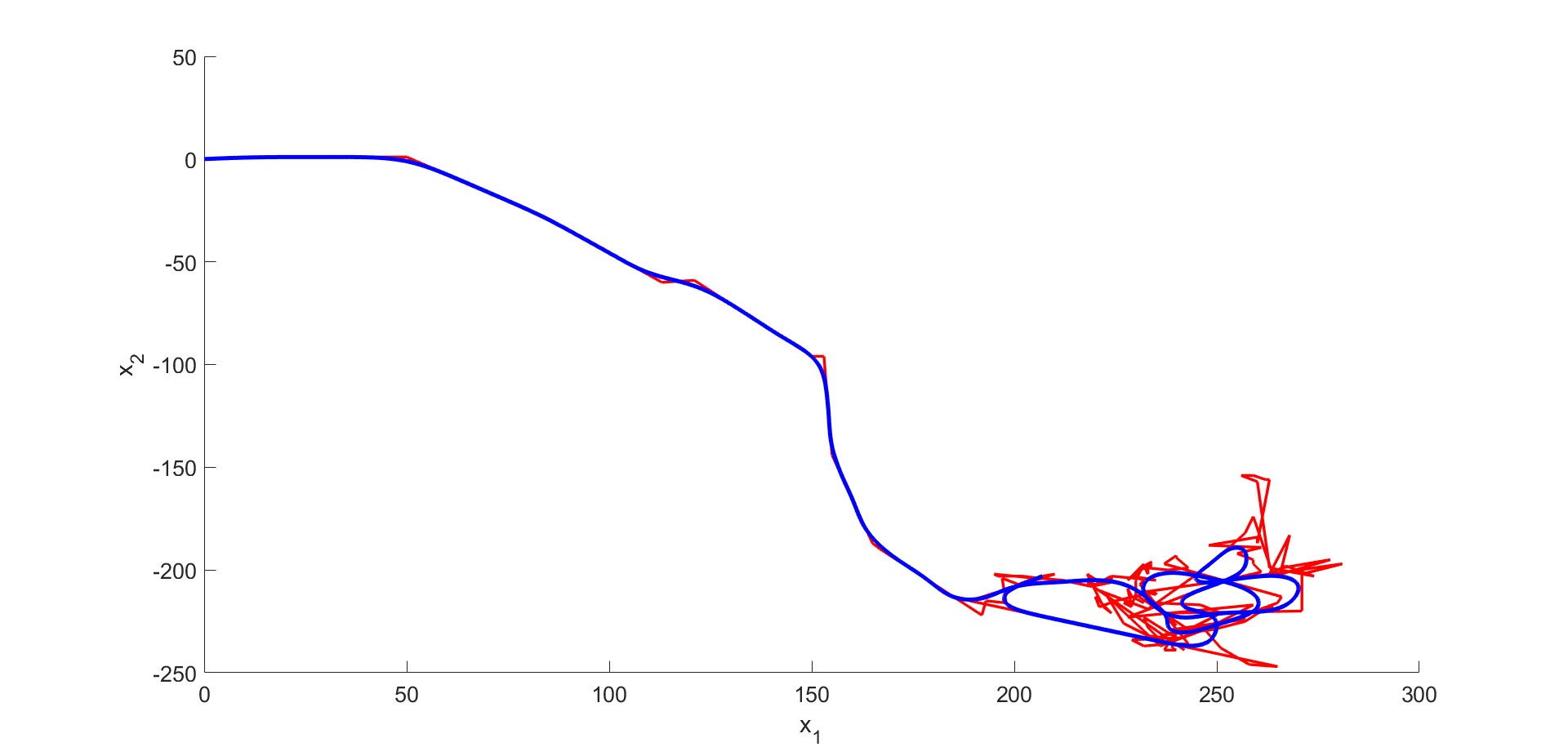}
  \qquad
  \caption{Original trajectory (red line) and smoothed trajectory (blue line). Results (from top to bottom) after $50$, $200$, $400$ and $580$ time steps.}
  \label{Fig5}
\end{figure}

\begin{figure}
  \includegraphics[width=10cm]{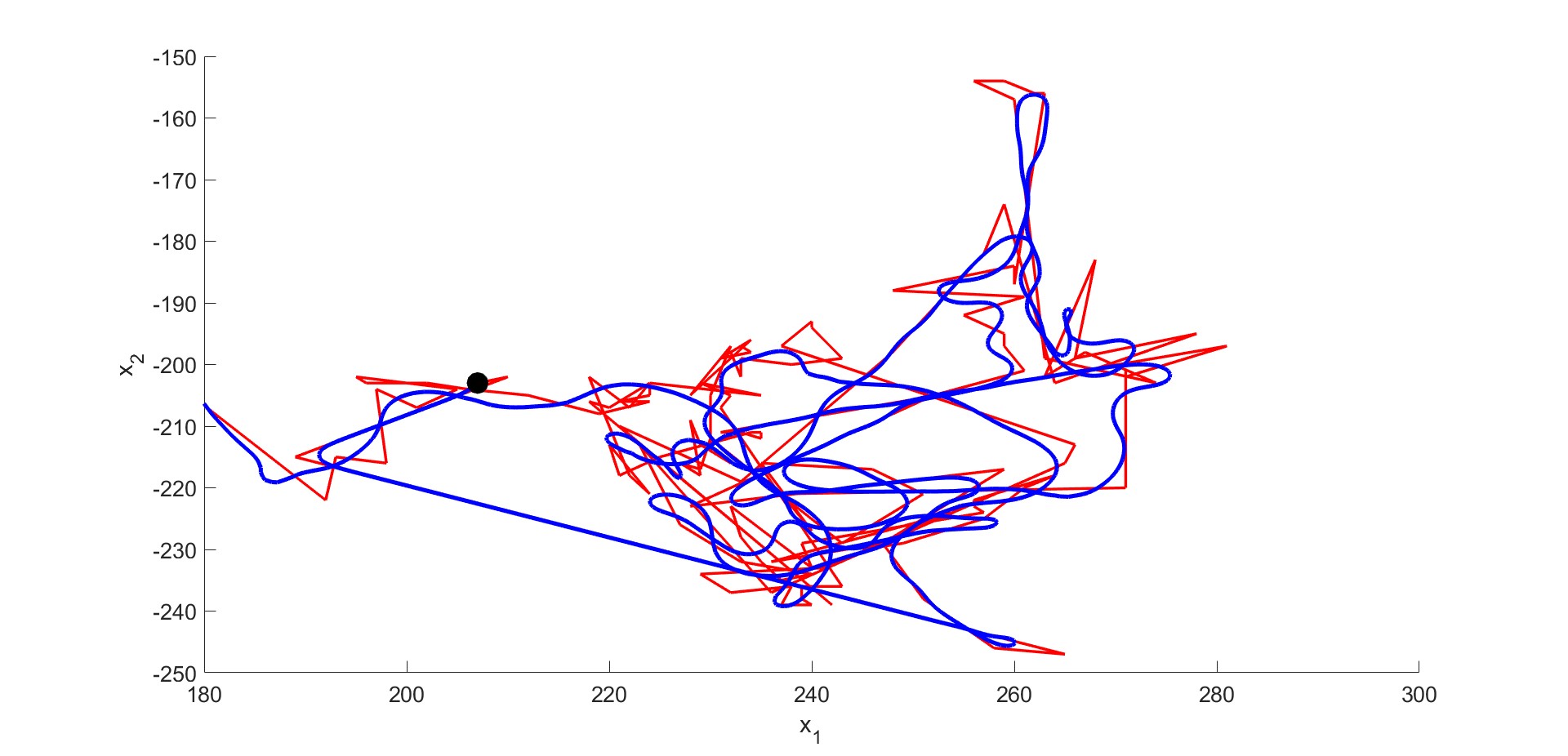}
  \includegraphics[width=10cm]{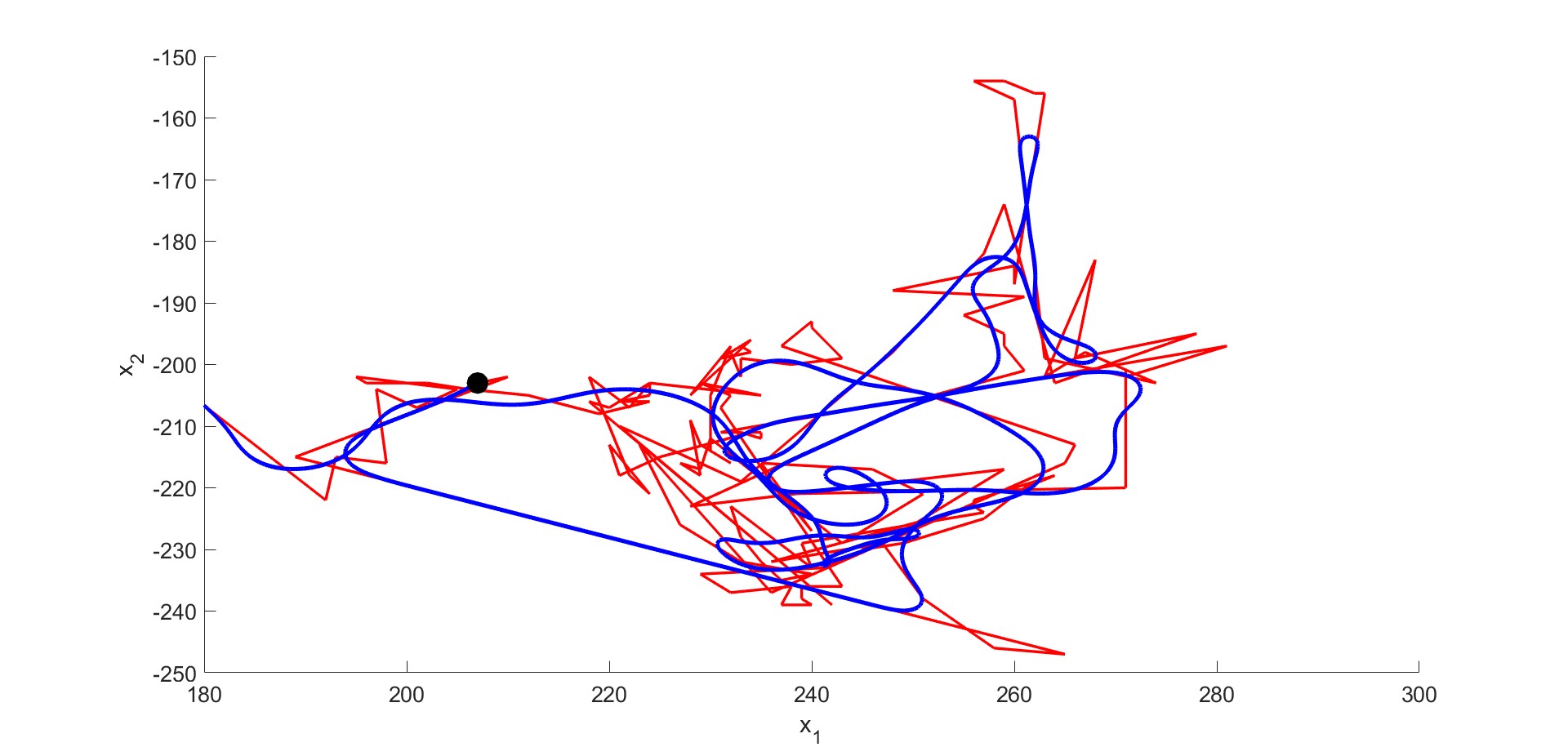}
  \includegraphics[width=10cm]{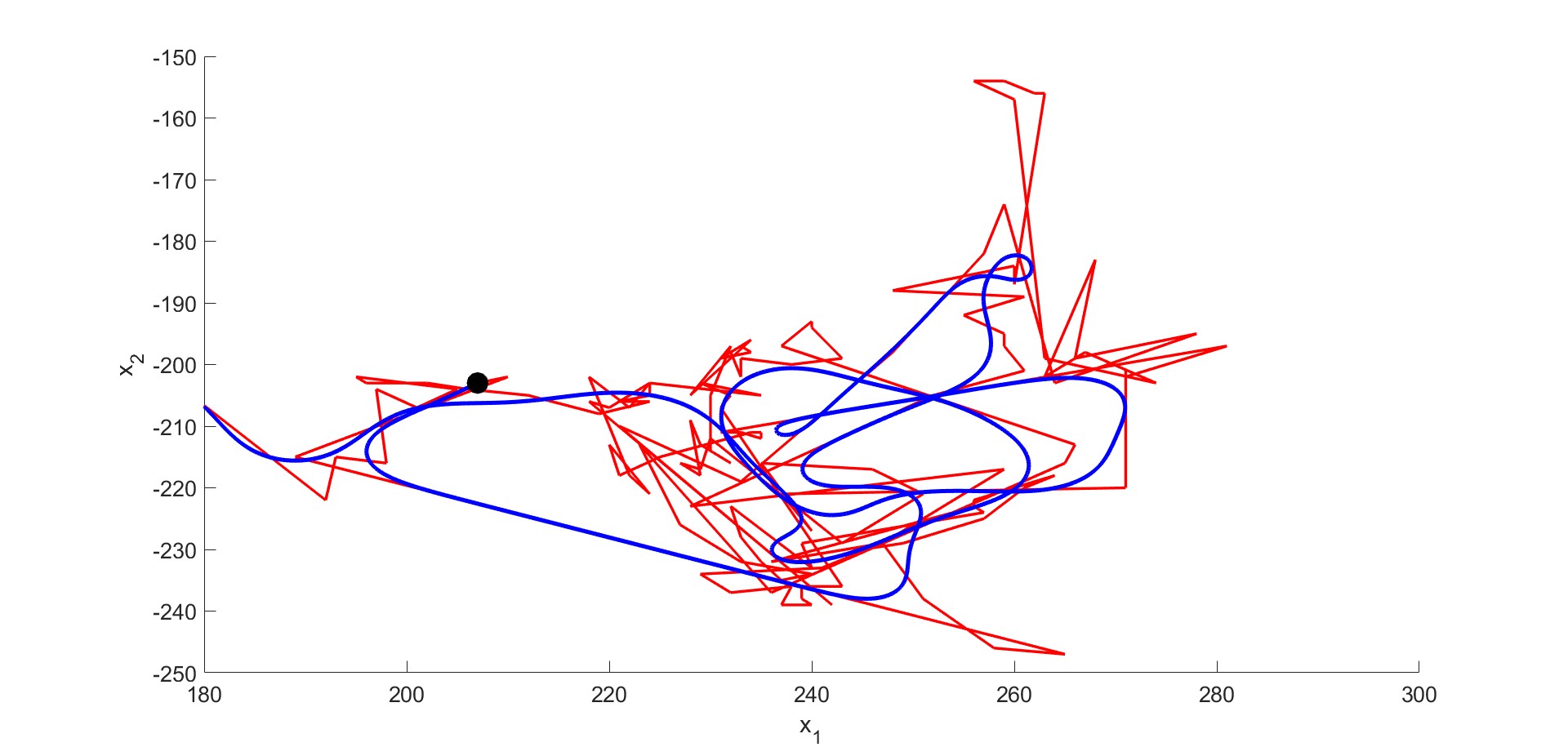}
  \includegraphics[width=10cm]{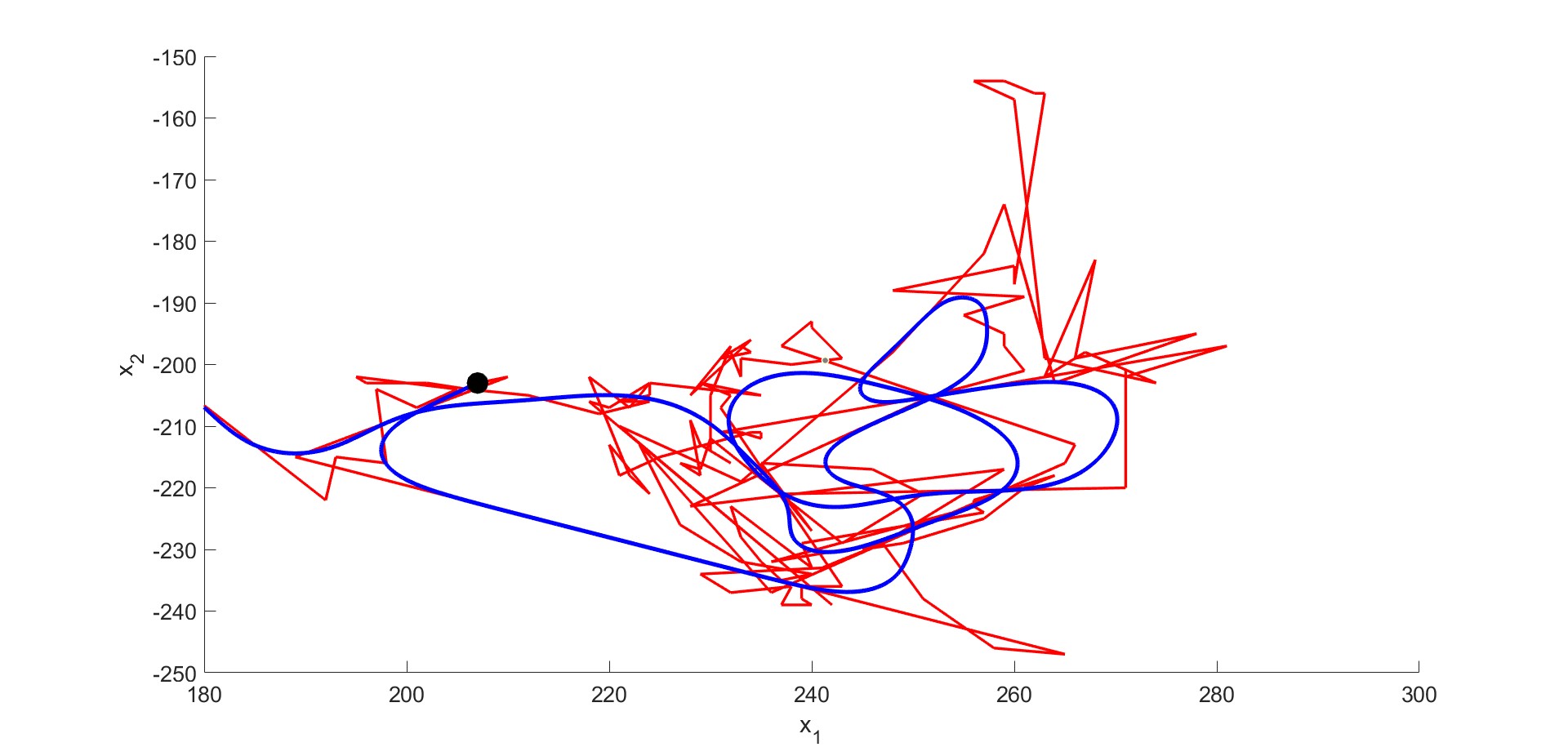}
  \qquad
  \caption{Detail of noisy part of trajectory in Fig.  \ref{Fig5}. Original trajectory (red line) and smoothed trajectory (blue line). Results (from top to bottom) after $50$, $200$, $400$ and $580$ time steps.}
  \label{Fig5zoom}
\end{figure}

\begin{figure}
  \includegraphics[width=10cm]{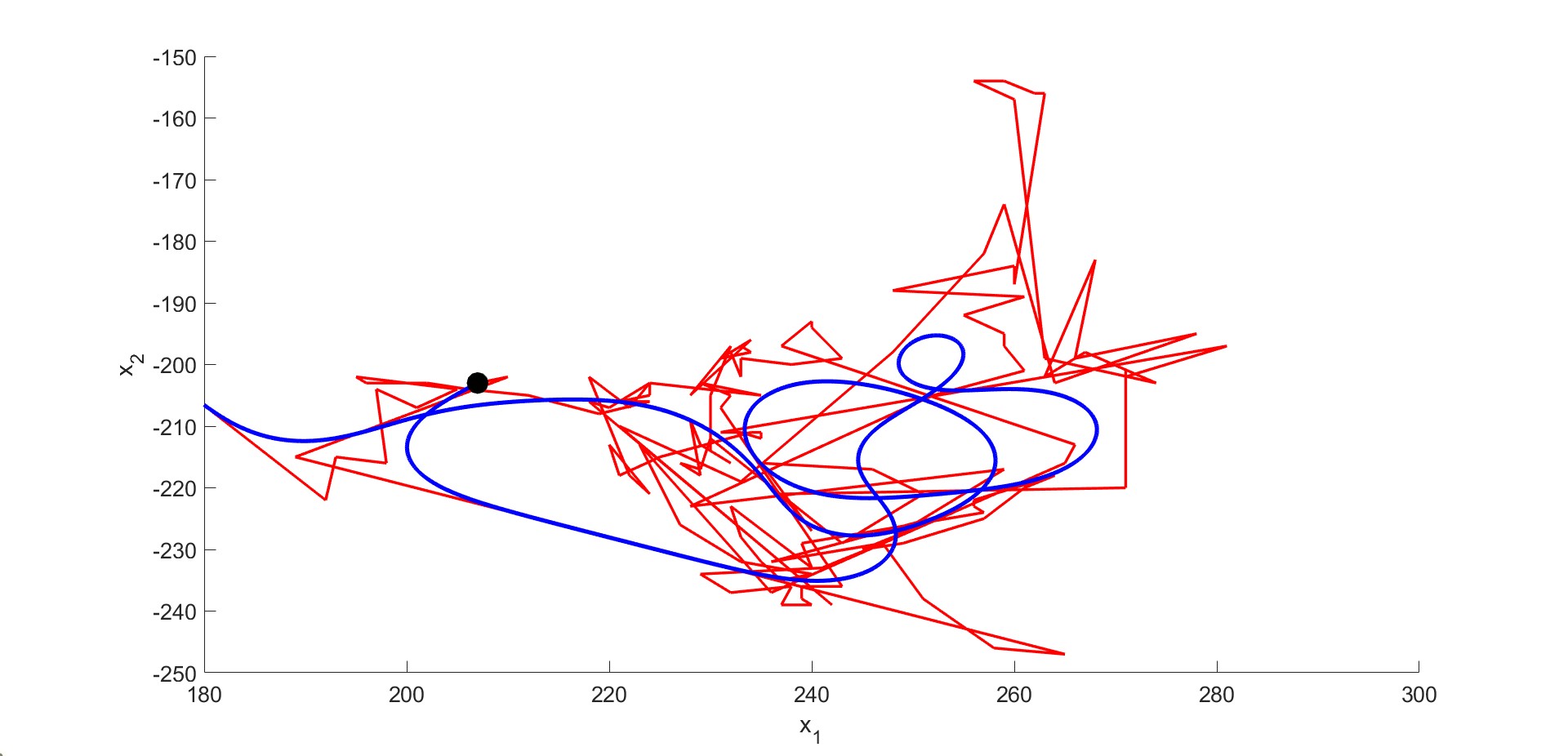}
  \includegraphics[width=10cm]{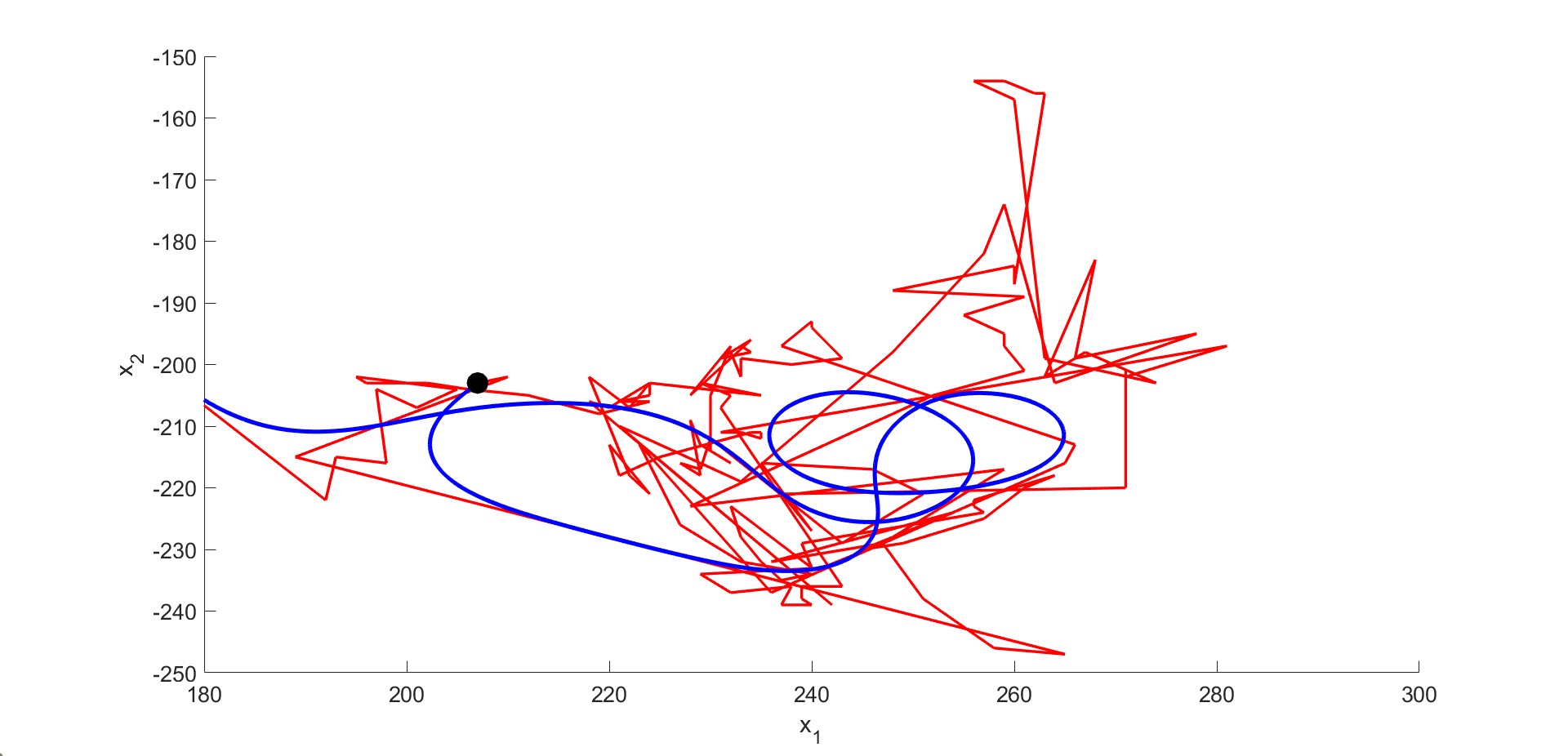}
  \includegraphics[width=10cm]{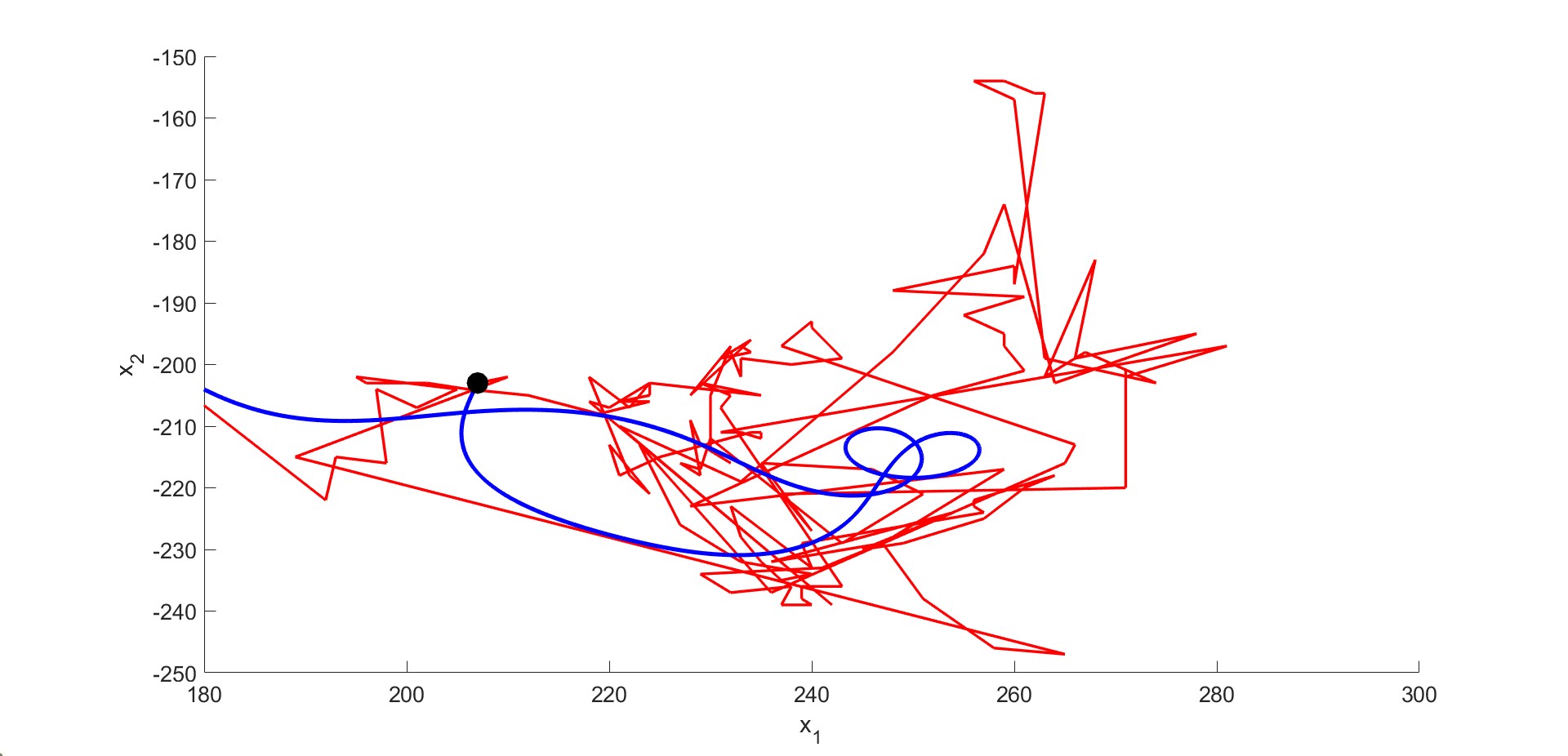}
  \includegraphics[width=10cm]{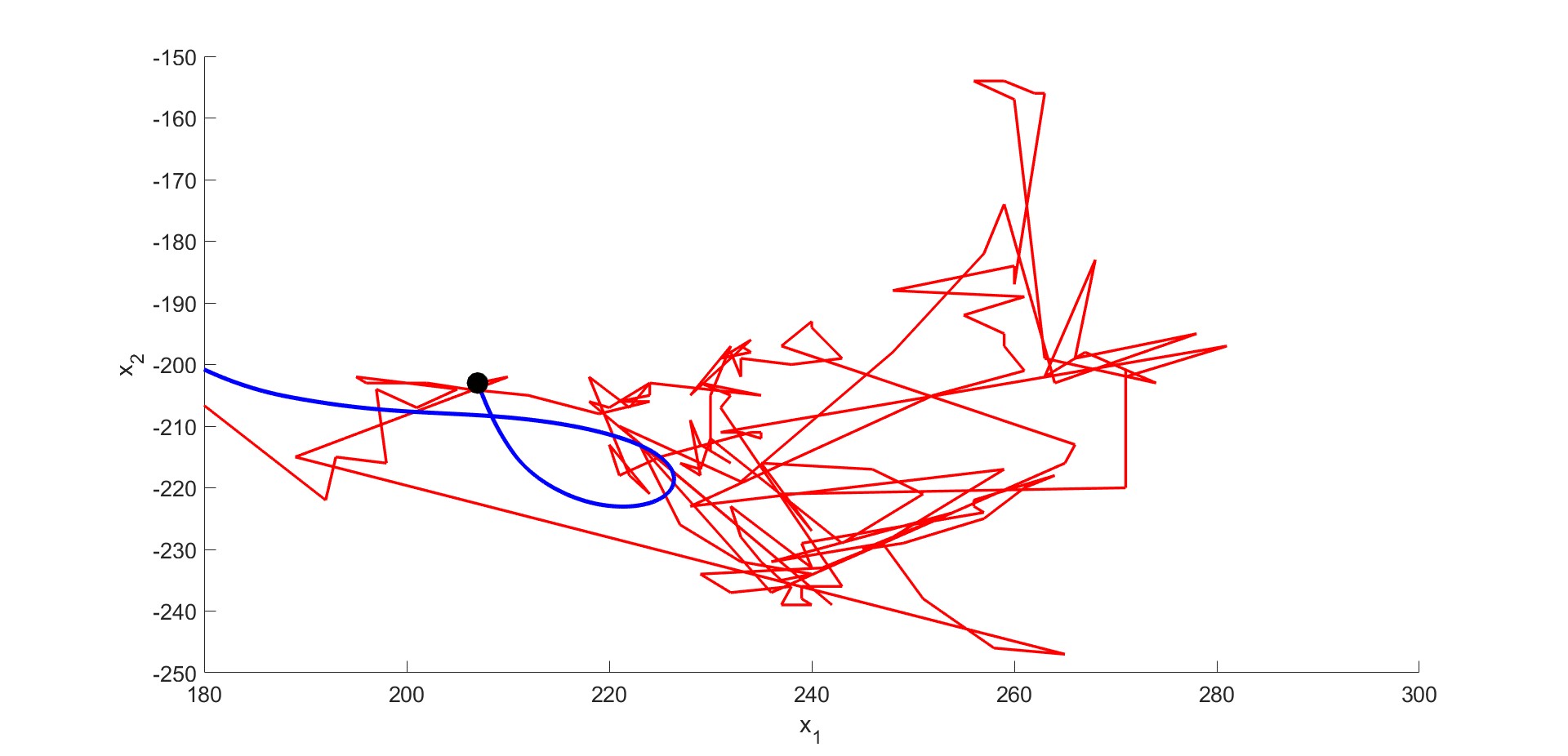}
  \qquad
  \caption{Detail of noisy part of trajectory in Fig.  \ref{Fig5}. Original trajectory (red line) and smoothed trajectory (blue line). Results (from top to bottom) after $1000$, $1500$, $2500$, and $5000$ time steps.}
  \label{Fig5zoom2}
\end{figure}

\begin{figure}
  \includegraphics[width=10cm]{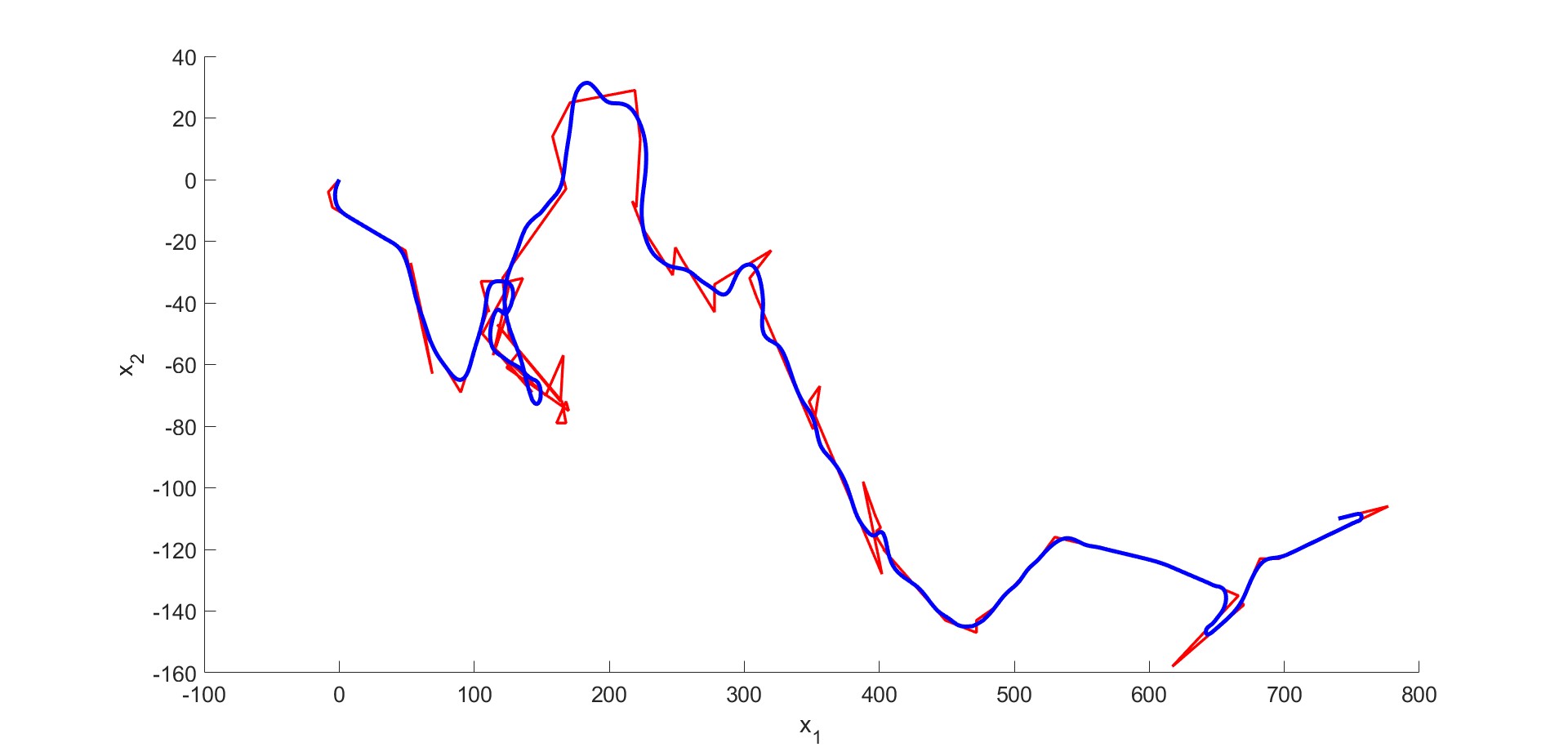}
  \includegraphics[width=10cm]{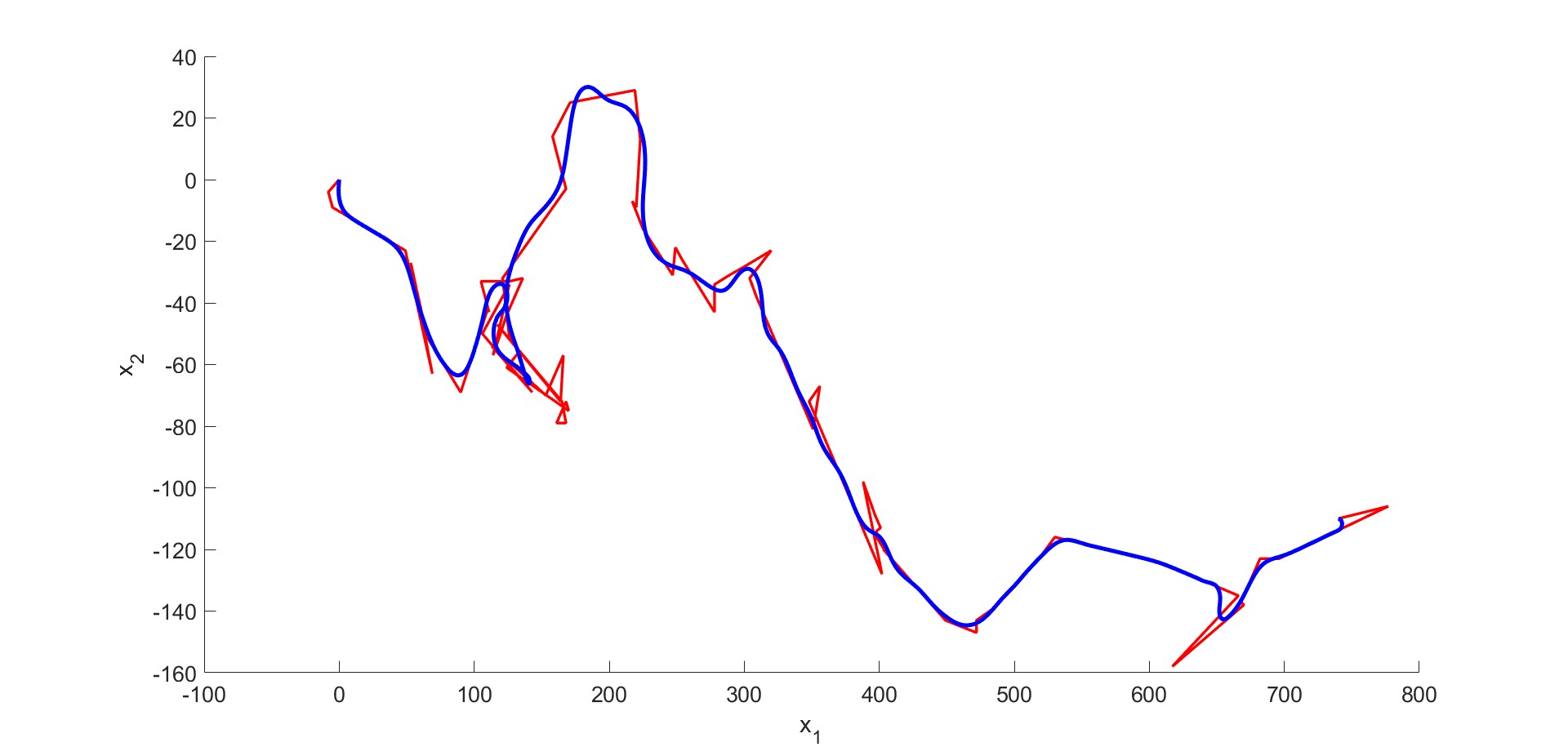}
  \includegraphics[width=10cm]{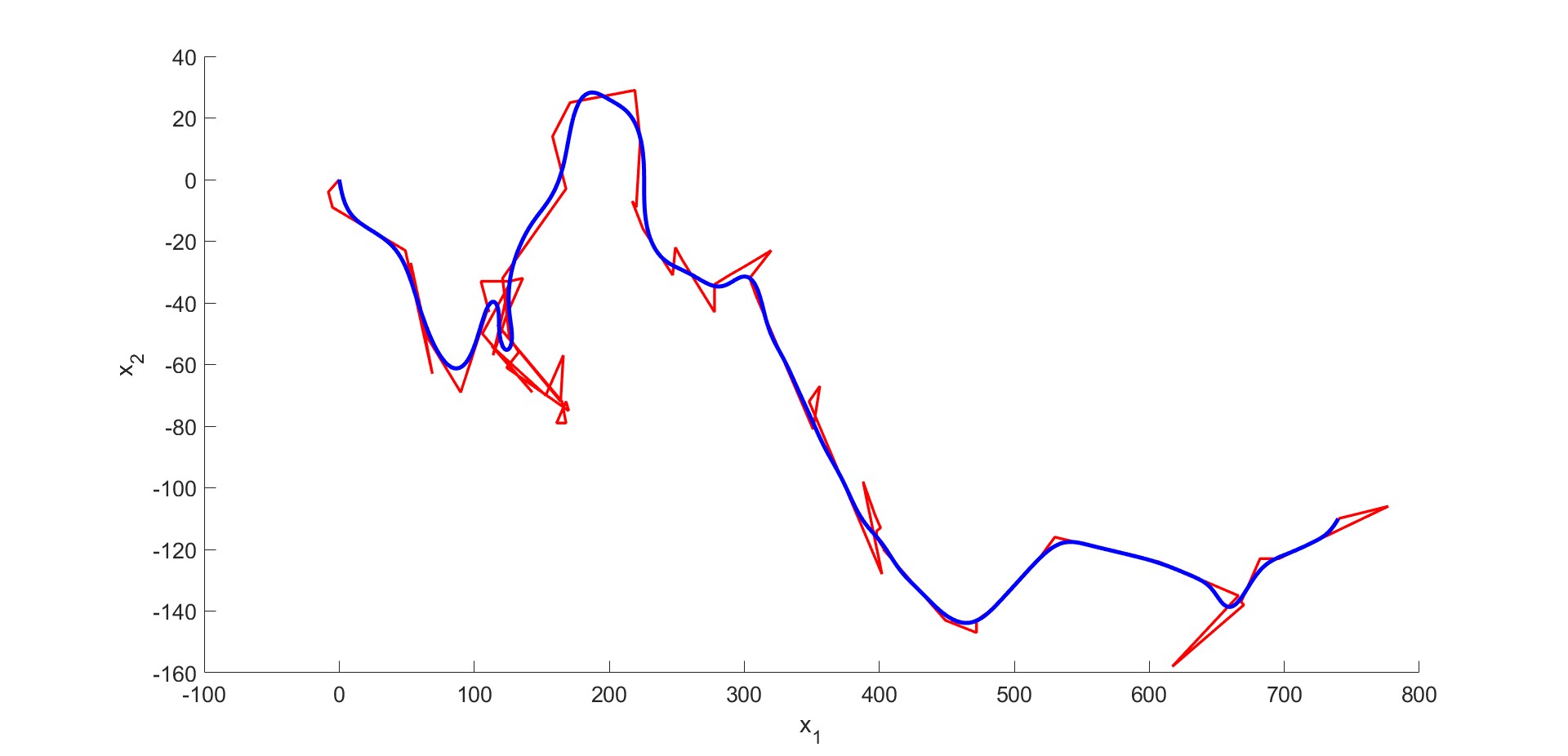}
  \includegraphics[width=10cm]{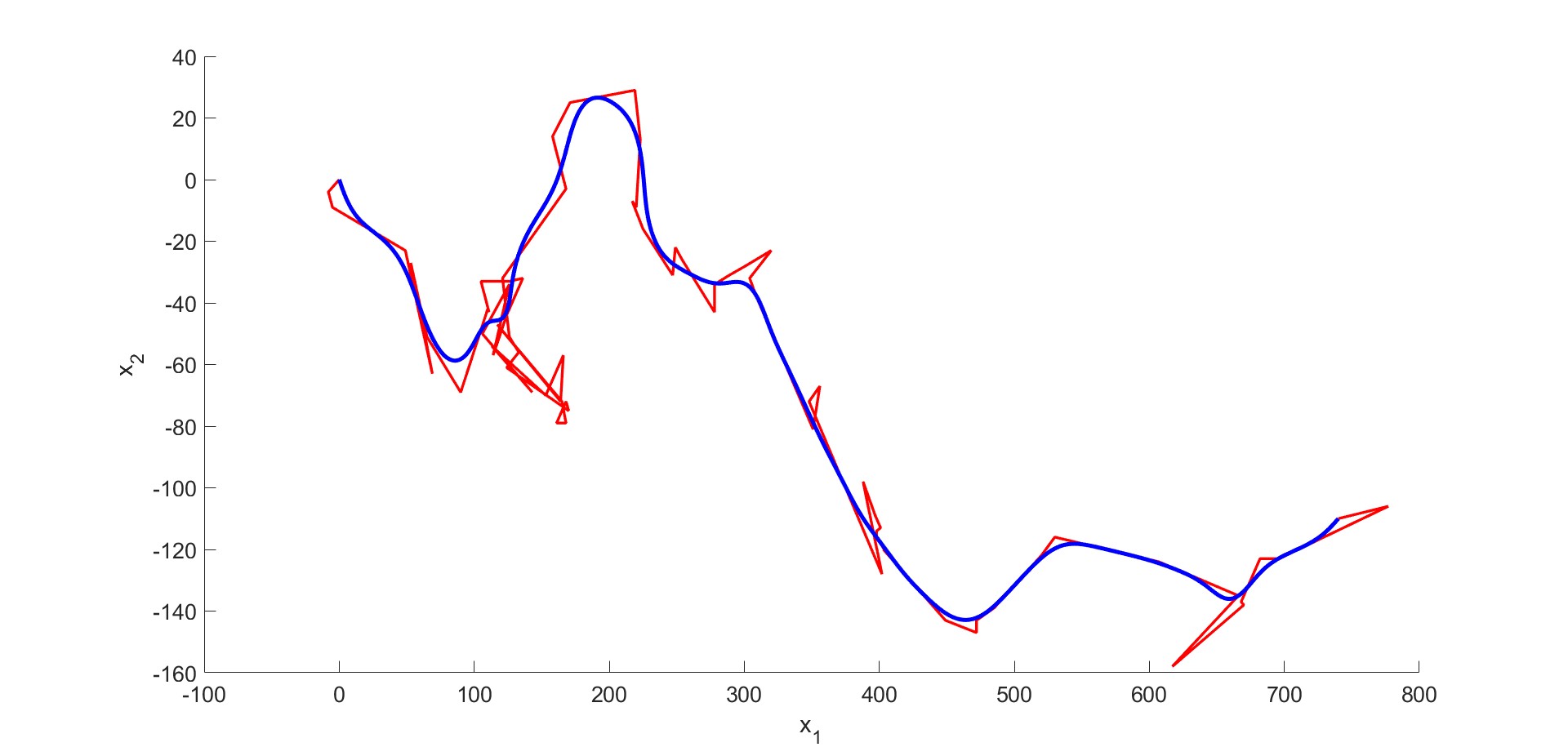}
  \qquad
  \caption{Original trajectory (red line) and smoothed trajectory (blue line). Results (from top to bottom) after $50$, $100$, $200$ and $340$ time steps.}
  \label{Fig6}
\end{figure}

\begin{figure}
  \includegraphics[width=10cm]{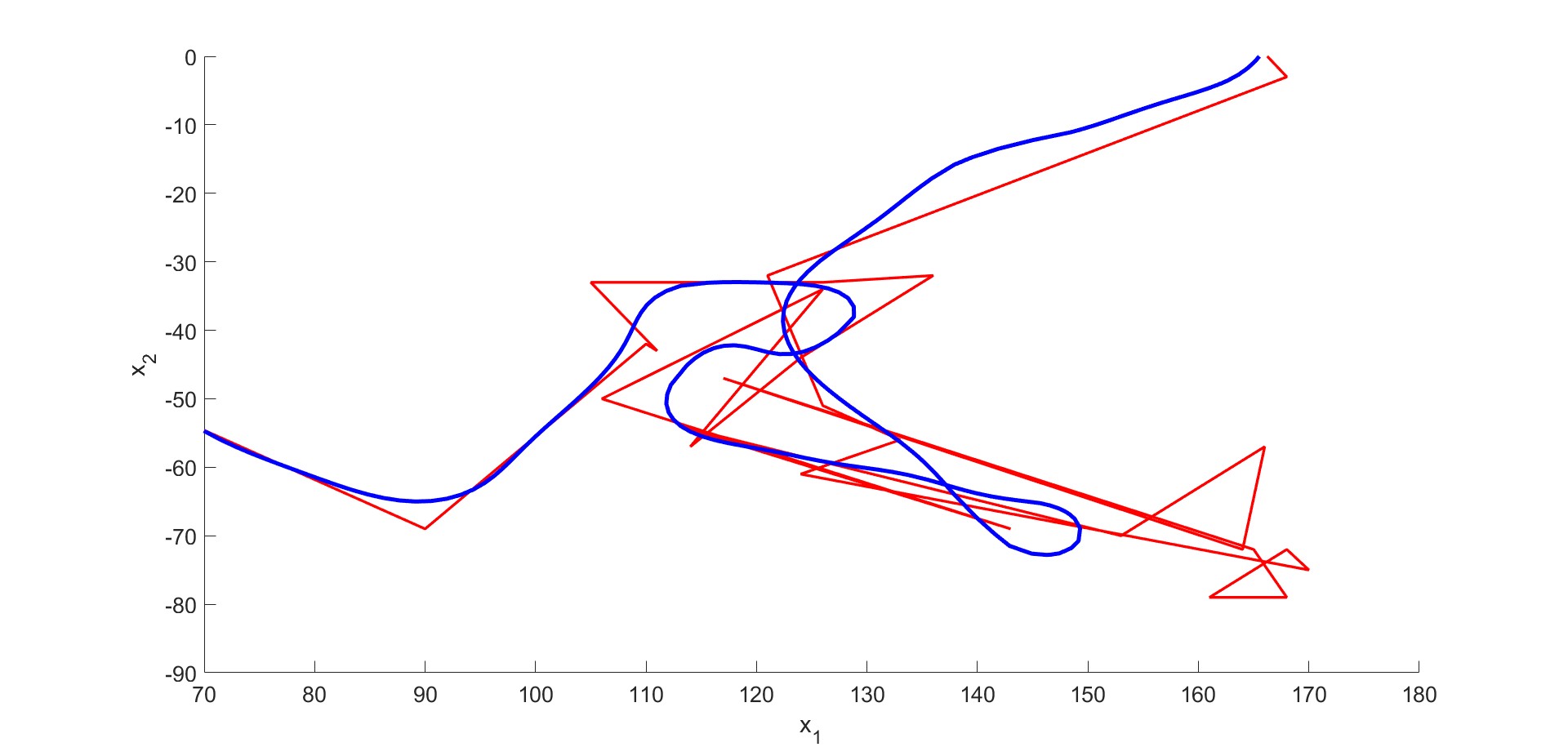}
  \includegraphics[width=10cm]{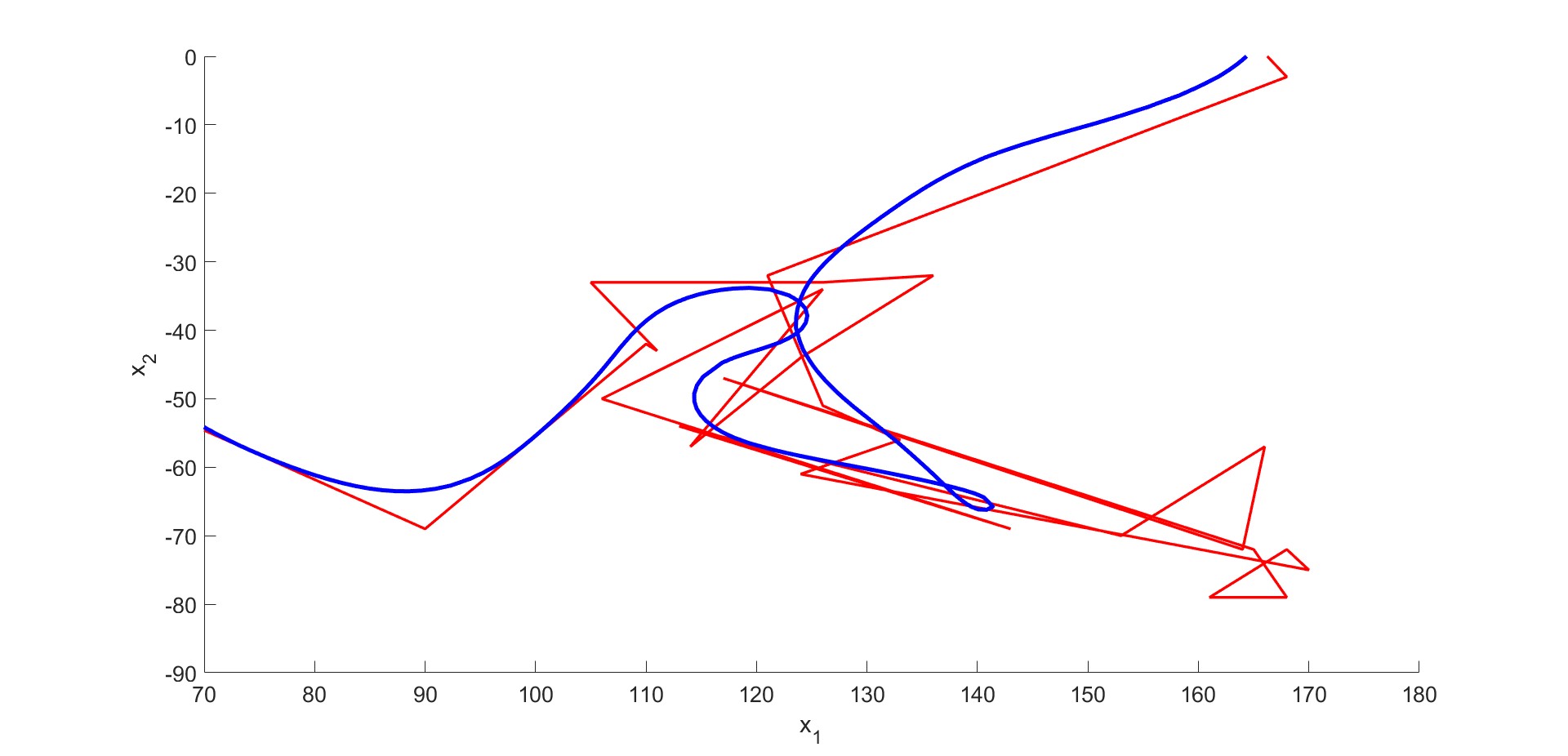}
  \includegraphics[width=10cm]{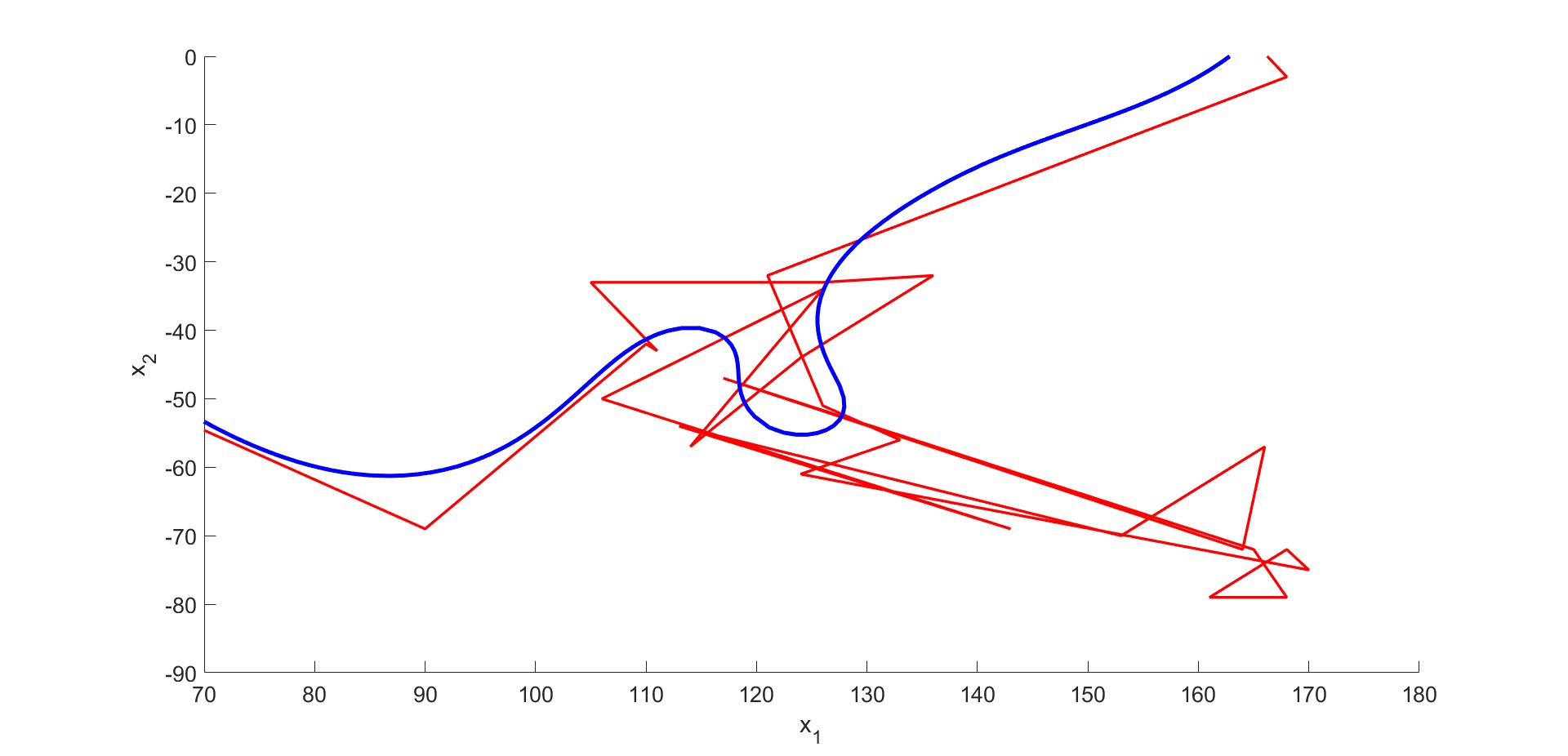}
  \includegraphics[width=10cm]{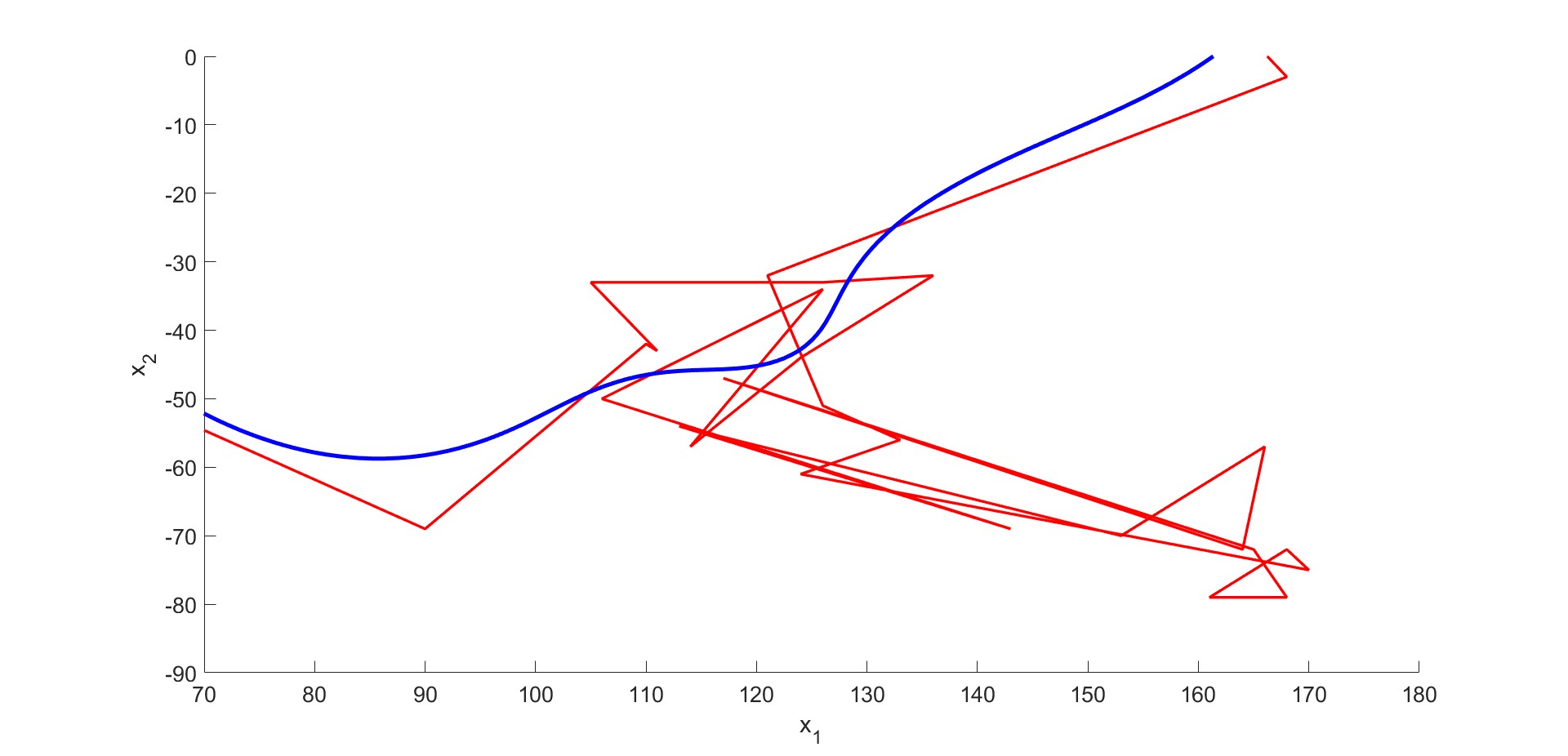}
  \qquad
  \caption{Detail of noisy part of trajectory in Fig.  \ref{Fig6}. Original trajectory (red line) and smoothed trajectory (blue line). Results (from top to bottom) after $50$, $100$, $200$ and $340$ time steps.}
  \label{Fig6zoom}
\end{figure}

Then, we focused on real data of macrophage motion. 
In the following experiments, we applied the proposed model to macrophage trajectories. The piecewise linear initial condition was found connecting the cell center of the macrophage for every frame of the $2$D$+$time data \cite{Sora}. First, we added grid points in the original trajectory so that we had almost uniformly distributed points. For every trajectory, the parameters were chosen as follows: $\lambda=1$, $\delta=0.005$, $\omega=1$, and $\tau=0.0001$. Figs. \ref{Fig5}-\ref{Fig6zoom} show the result of the smoothing for two different trajectories. 
The red lines in Figs. \ref{Fig5} and \ref{Fig6} are the original trajectories for two different macrophages. As one can observe, due to the random motion of the cells, these trajectories usually have noisy parts: Figs \ref{Fig5zoom}, \ref{Fig5zoom2}, and \ref{Fig6zoom} show the details of these parts. The blue lines are the smoothed trajectories after different numbers of time steps: the evolving curve is kept close to the original trajectory by the attracting term but we also observe the smoothing of the noisy parts due to the diffusion term. In Figs. \ref{Fig5zoom} and \ref{Fig5zoom2}, we show the result of the smoothing until the random part of the trajectory is completely smoothed: the black dot represents the last point of the trajectory. \\\\  
As a stopping criterion, we considered the mean Hausdorff distance between two discrete curves. Consider the discrete sets of points $\mathcal{A}=\{a_0,...,a_{n+1}\}$, $\mathcal{B}=\{b_0,...,b_{n+1}\}$: we will indicate the elements of the discrete curves by $\textbf{a}_i=a_i-a_{i-1}$ and  $\textbf{b}_i=b_i-b_{i-1}$, $i=1,...,n+1$ respectively. Let us indicate by $A=\{\textbf{a}_1,...,\textbf{a}_{n+1}\}$ and $B=\{\textbf{b}_1,...,\textbf{b}_{n+1}\}$ the sets of elements. The mean Hausdorff distance $\bar{d}_H(\mathcal{A}, \mathcal{B})$ is then defined as
\begin{equation}
\begin{split}
\bar{\delta}_H(\mathcal{A}, \mathcal{B})&=\frac{1}{n}\sum_{i=1}^{n}\min_{\textbf{b}\in B}d(a_i,\textbf{b}),\\
\bar{\delta}_H(\mathcal{B}, \mathcal{A})&=\frac{1}{n}\sum_{i=1}^{n}\min_{\textbf{a}\in A}d(b_i,\textbf{a}),\\
\bar{d}_H(\mathcal{A},\mathcal{B})&=\frac{\bar{\delta}_H(\mathcal{A}, \mathcal{B})+\bar{\delta}_H(\mathcal{B}, \mathcal{A})}{2},\\
\end{split}
 \label{HD}
\end{equation}
where $d(a_i,\textbf{b})$ is the distance between a point $a_i \in \mathcal{A}$ and an element $\textbf{b}\in B$ defined in section \ref{S2.2}. Since in our model the endpoints $a_0,~a_{n+1},~b_0,~b_{n+1}$ are fixed, we do not consider them in the computation of the mean Hausdorff distance. In order not to slow down the computation too much this distance was calculated only every $p$ number of time steps; we chose the tolerance $\epsilon=0.000065$ and stopped iterating when $d_H(\textbf{x}_n,\textbf{x}_{n+p})<\epsilon$. 
For the trajectory in Fig. \ref{Fig5}, the number of time steps needed to stop the evolution was $580$, while for the one in Fig. \ref{Fig6}, it was $340$.
\subsection{Velocities on smoothed curves}
\label{S3.2}
\begin{figure}
  \centering
     \includegraphics[width=12cm,center]{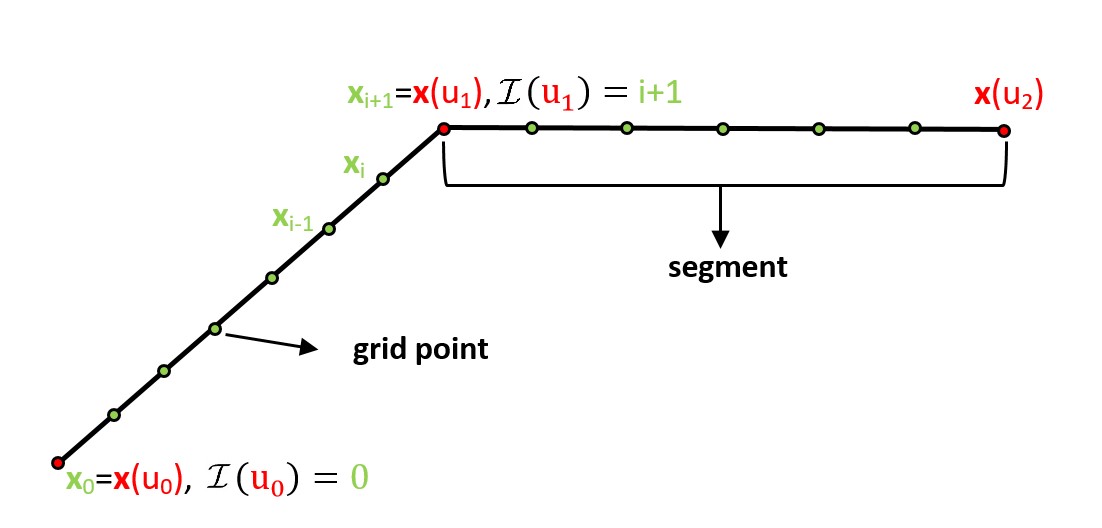}
  \caption{Visualization of curve discretization. The green dots represent the grid points and the red dots represent the endpoints of the segments.}
  \label{I_expl}
\end{figure}

Our final goal was to compute the velocity for every grid point on the evolved curve. In the $2$D$+$time data, the real-time between each frame was $1$ time unit between the frames, which is equal to 2 minutes. That means, in the original piecewise linear trajectory, the time between two endpoints of a segment was $1$ time unit (see Fig. \ref{I_expl}). First, when we added the grid points inside the initial segments, we did it in such a way that the new added points were uniformly distributed inside the segment. Consequently, we calculated the original velocity for the grid points as the original length of the segment divided by the time interval $\Delta_t$, that in our case was $1$ (time unit), and consider it constant inside the segment. Therefore, in the beginning, all the added grid points inside a segment have the same velocity. We assume this property holds also for the new velocities of the grid points on the evolved curve. Our goal was to derive the velocity for every grid point on the evolved curve but we couldn't do it directly because, due to the stability of the numerical computations, we added the tangential velocity in our model and this caused an asymptotically uniform redistribution of the points. If the trajectory is not too complicated, the tangential velocity does not move the points too much. But, in our case, some of the trajectories had noisy parts with many points, which were moved a lot during the evolution, so we were losing almost completely the spatial information.\\
Then, to get the new velocity estimation, first we studied the evolution of the length of the segments. For the total length, we have already shown in section \ref{S2.1} that it holds
\begin{equation}
 L_t=\int_{0}^{1}k\beta ds+\int_{0}^{1} \alpha_s ds=\int_{0}^{1}k\beta ds+\alpha(1)-\alpha(0).
\end{equation}
In the case of the total length, we had fixed the endpoints of the trajectory so we had $\alpha(1)=\alpha(0)=0$. Similarly, for the evolution of every segment, we can consider the same formula. We will denote by $j$ the index for the segments and by $i$ the index for the grid points. Consider the $j$-th segment: we will indicate by $u_{j-1}$ and $u_{j}$ the first and last endpoint of the segment parametrization, respectively. Therefore, for the evolution of the length $L_j$ of the $j$-th segment, it holds
\begin{equation}
 (L_j)_t=\int_{u_{j-1}}^{u_j}k\beta ds+\int_{u_{j-1}}^{u_j} \alpha_s ds=\int_{u_{j-1}}^{u_j}k\beta ds+\alpha(u_{j})-\alpha(u_{j-1}).
\end{equation}
Notice that in this case, in general, we have $\alpha(u_j)\neq 0$ and $\alpha(u_{j-1})\neq0$. But, as we have already underlined, the tangential velocity is used only for the stability of the numerical computations and it does not influence the shape of the evolving curve. So, to not lose the spatial information, we will calculate the new length of the segment as if we have applied to the trajectory the model with $\alpha=0$, i.e., $\alpha(u_j)=\alpha(u_{j-1})=0$. Consequently, we consider the following model for the evolution of the $j$-th segment 
\begin{equation}
 (L_j)_t=\int_{u_{j-1}}^{u_j}k\beta ds.
 \label{EQ3.34}
\end{equation}
Let's now consider the discretization: let $m$ be the time step index and $\tau$ the length of the discrete time step. We approximate the time derivative by the finite difference and consider the formulas for $k,~\beta$ defined in (\ref{EQ2.30}). To relate the index $j$ of the segments to the index $i$ of the grid points we will indicate by $\mathcal{I}(u_{j-1})$, $\mathcal{I}(u_j)$ the corresponding $i$ index of the endpoint $\textbf{x}(u_{j-1})$, respectively $\textbf{x}(u_j)$ (see Fig. \ref{I_expl}). Then, we obtain
\begin{equation}
 \frac{L_{j}^{m+1}-L_{j}^{m}}{\tau}= \sum_{i=\mathcal{I}(u_{j-1})+1}^{\mathcal{I}(u_j)}h_{i}^{m}k_{i}^{m}\beta_{i}^{m}.
\end{equation}
Notice that we consider the $i$ index starting from $\mathcal{I}(u_{j-1})+1$: that is because $h_{\mathcal{I}(u_{j-1})+1}$ is defined as $h_{\mathcal{I}(u_{j-1})+1}=\lvert \textbf{x}_{\mathcal{I}(u_{j-1})+1}-\textbf{x}_{\mathcal{I}(u_{j-1})}\rvert$, which is indeed the first element $h_i$ of the segment with endpoints $\textbf{x}(u_{j-1})$ and $\textbf{x}(u_{j})$. We obtain the fully discrete approximation of equation (\ref{EQ3.34}) for the new length of the $j$-th segment
\begin{equation}
L_{j}^{m+1}=L_{j}^{m}+\tau \sum_{i=\mathcal{I}(u_{j-1})+1}^{\mathcal{I}(u_j)}h_{i}^{m}k_{i}^{m}\beta_{i}^{m}.
\label{EQ3.36}
\end{equation}
For the first and last element of each segment, we use a different computation of $k$ and $\beta$ that does not consider informations from the previous and next segment.\\
If we consider the evolution of the curve driven only by the curvature influence, then $\beta_i=-\delta k_i$ and we obtain
\[
 \sum_{i=\mathcal{I}(u_{j-1})+1}^{\mathcal{I}(u_j)}h_{i}^{m}k_{i}^{m}\beta_{i}^{m}=-\sum_{i=\mathcal{I}(u_{j-1})+1}^{\mathcal{I}(u_j)}h_{i}^{m}\delta(k_{i}^{m})^2<0.
\]
Therefore $L_{j}^{m+1}<L_{j}^{m}$ for every segment $j$ and for every time step index $m$. Experimentally, this property holds also if $\beta_i=-\delta k_i+\lambda w_i$. If the segment is in a part of the trajectory which has high curvature (for example in the noisy parts where the cell shows random motion),  the length of that segment will decrease faster than the length of a segment in a part of the trajectory with small curvature (for example in the parts where the cell shows directional motion). It is then clear that some of the segments may eventually disappear and we may obtain $L_{j}^{m+1}=0$.\\\\
Let us first focus on the simple case where the length of one segment is decreasing but the segment is not disappearing. In the first time step, $\mathcal{I}(u_{j-1})$ and $\mathcal{I}(u_{j})$ are loaded from the original trajectory and are the indexes of the endpoints of the original segments. Later, the length of every segment changes so, after applying the formula in (\ref{EQ3.36}), we need to find the new $i$ indexes for $\mathcal{I}(u_{j-1})$ and $\mathcal{I}(u_j)$. If we indicate by
\begin{equation}
 \begin{split}
  L_{s}^{m+1}&=\sum_{j=1}^{M-1} L_{j}^{m+1}\\
 \end{split}
\end{equation}
where $M-1$ is the number of original segments. Ideally, $L_{s}^{m+1}= L^{m+1}$, where $L^{m+1}$ is defined in (\ref{EQ2.30}).  However, in the numerical computations, in later time steps, it holds $L_{s}^{m+1}\neq L^{m+1}$. We define the ratio of the $j$-th segment as
\begin{equation}
 r_{j}^{m+1}=\frac{L_{j}^{m+1}}{L_{s}^{m+1}},
\end{equation}
and the discrete length of the segment is defined as 
\begin{equation}
 L^{d}_{j}=r_{j}^{m+1}L^{m+1},
\end{equation}
so that $\sum_{j=1}^{M-1}L^{d}_{j}=L^{m+1}$. Once we obtain the new lengths $L^{d}_{j}$, to get the new $i$ indexes $\mathcal{I}(u_{j-1})$ and $\mathcal{I}(u_{j})$ for the endpoints of the segments we apply the following algorithm. Set $i=1$, $j=1$, $S=L^{d}_{j}$, $S_e=0$ and repeat until $j=M-1$ and $i=n+1$:
\begin{itemize}
 \item $S_e=S_e+h_{i}$\\
 \item if $S_e\geq S$
 \begin{itemize}
  \item if $S_e> S$ then $\mathcal{I}(u_{j})=i-1$\\
 \item if $S_e=S$ then $\mathcal{I}(u_{j})=i$\\
 \item update $S=S+L^{d}_{j+1}$.\\
 \item $j=j+1$\\
 \end{itemize}
 \item $i=i+1$\\
\end{itemize}
Notice that we focused only on the index of the last endpoint of each segment $\mathcal{I}(u_j)$: that is because the last endpoint of the segment $j$ is the first endpoint of the segment $j+1$.\\

\begin{figure}
  \includegraphics[width=12cm]{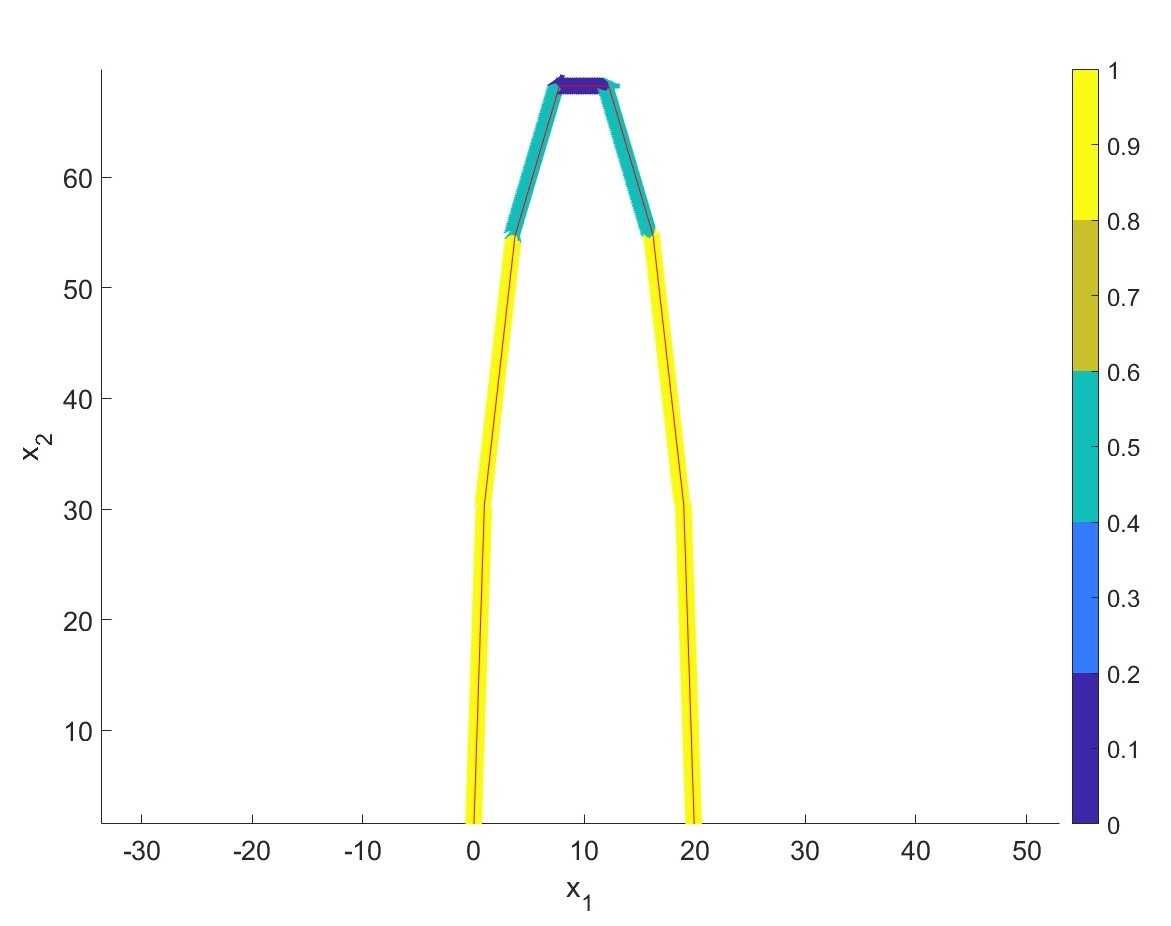}
  \includegraphics[width=12cm]{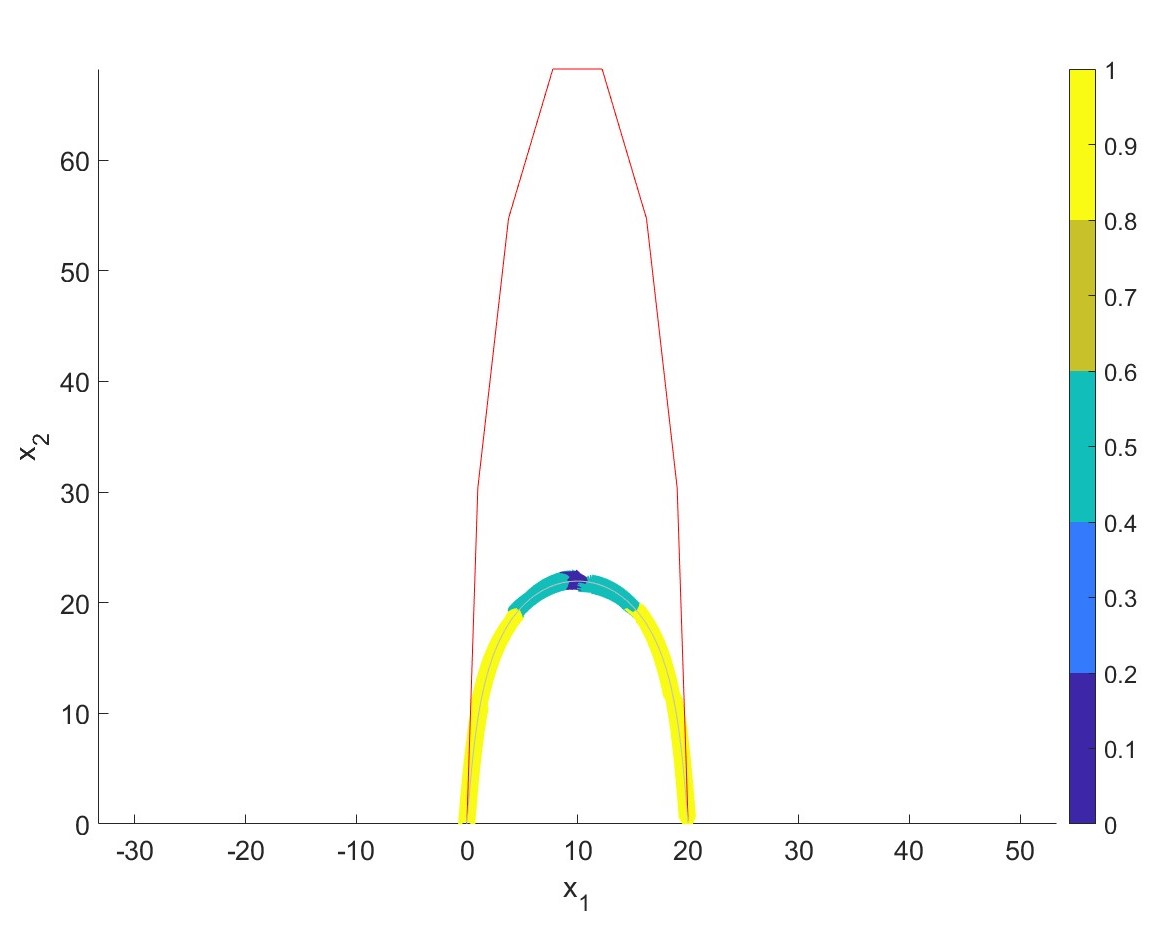}
  \qquad
  \caption{Top: initial trajectory and plotting of the velocity vectors. Bottom: initial trajectory (red line), evolved curve, and plotting of the velocity vectors. The vectors are normalized for better visualization, color represents the norm.}
  \label{Fig7}
\end{figure}

\begin{figure}
  \includegraphics[width=13cm]{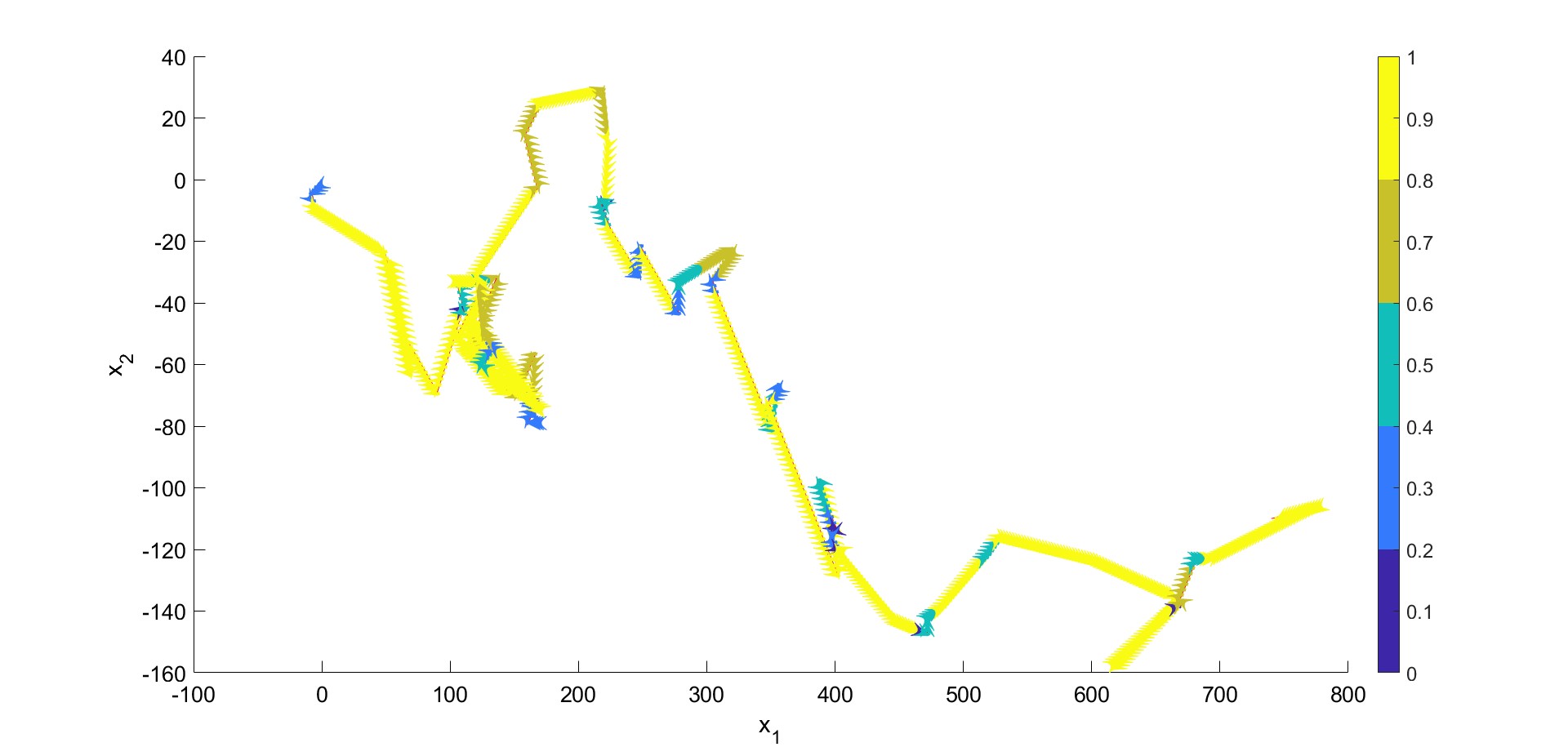}
  \includegraphics[width=13cm]{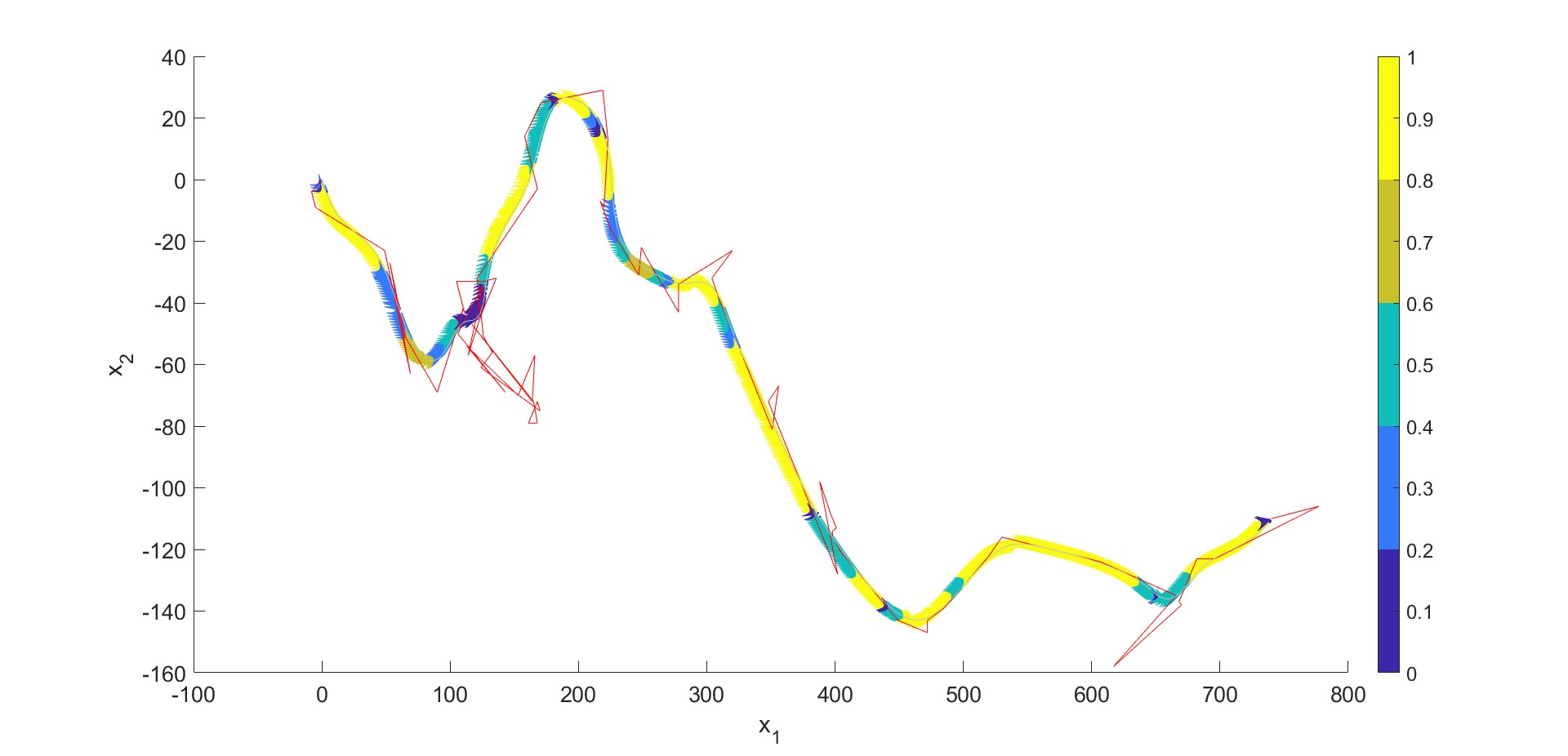}
  \qquad
  \caption{Top: initial trajectory and plotting of the velocity vectors. Bottom: initial trajectory (red line), evolved curve, and plotting of the velocity vectors. The vectors are normalized for better visualization, color represents the norm.}
  \label{Fig8}
\end{figure}

\begin{figure}
  \includegraphics[width=13cm]{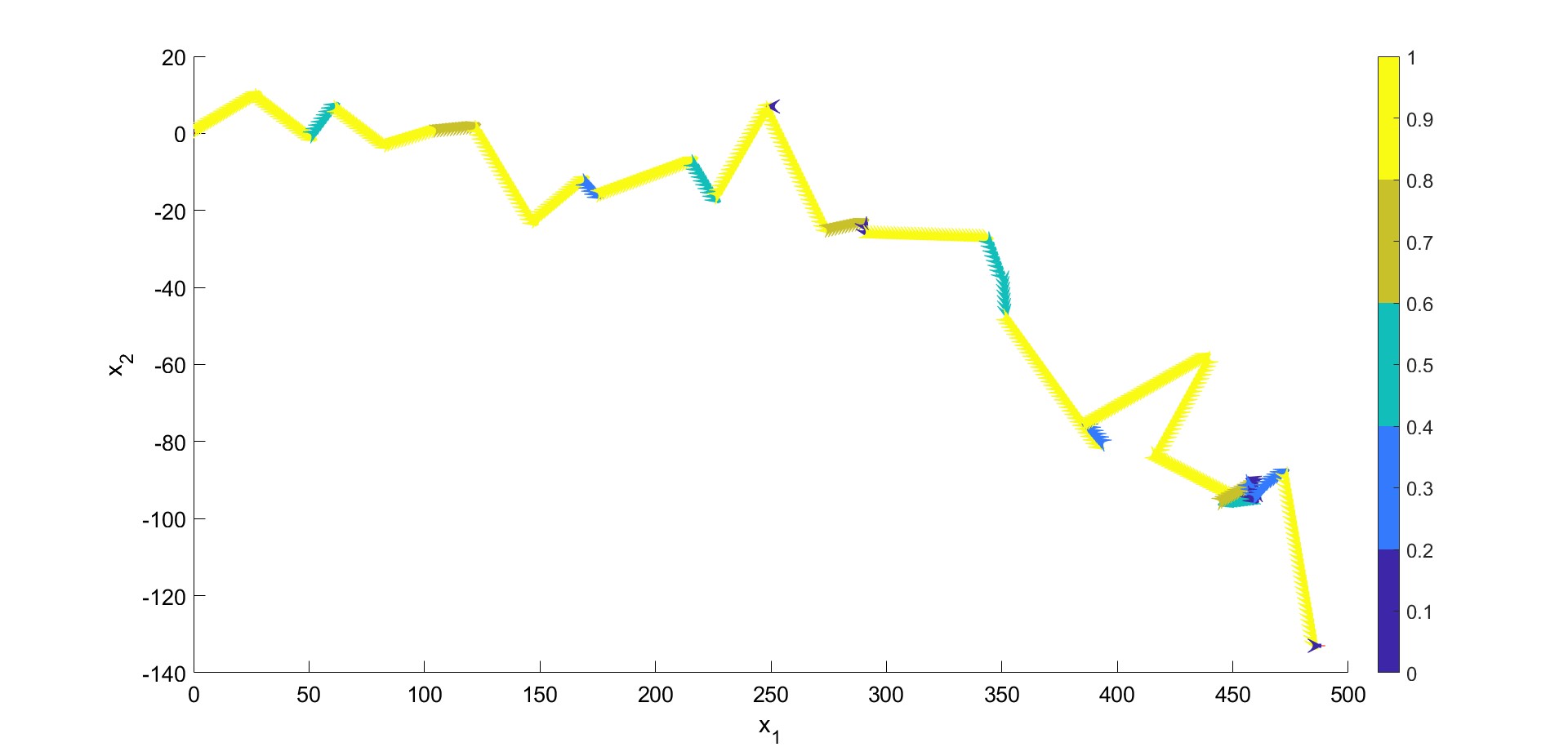}
  \includegraphics[width=13cm]{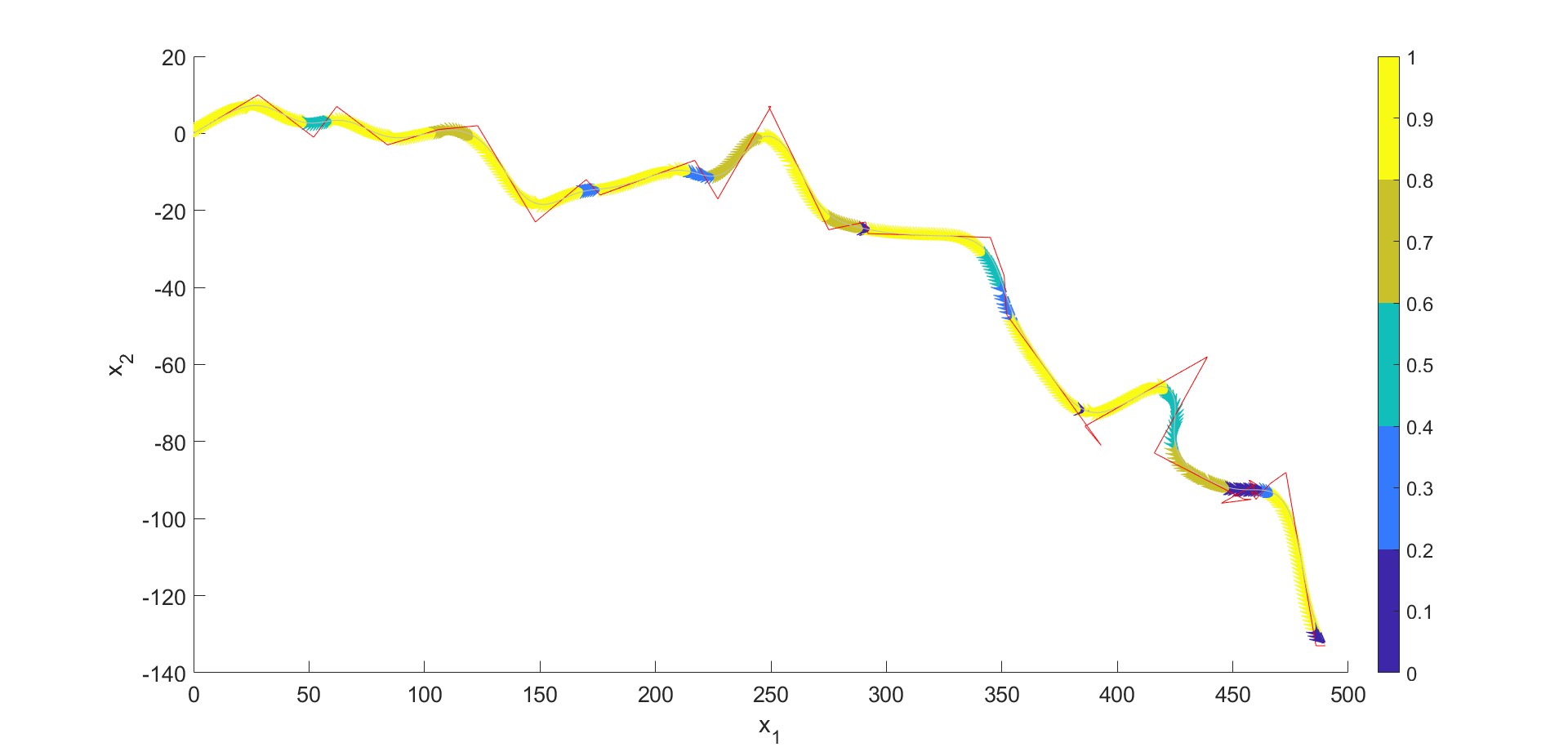}
  \qquad
  \caption{Top: initial trajectory and plotting of the velocity vectors. Bottom: initial trajectory (red line), evolved curve, and plotting of the velocity vectors. The vectors are normalized for better visualization, color represents the norm.}
  \label{Fig9}
\end{figure}

\begin{figure}
  \includegraphics[width=13cm]{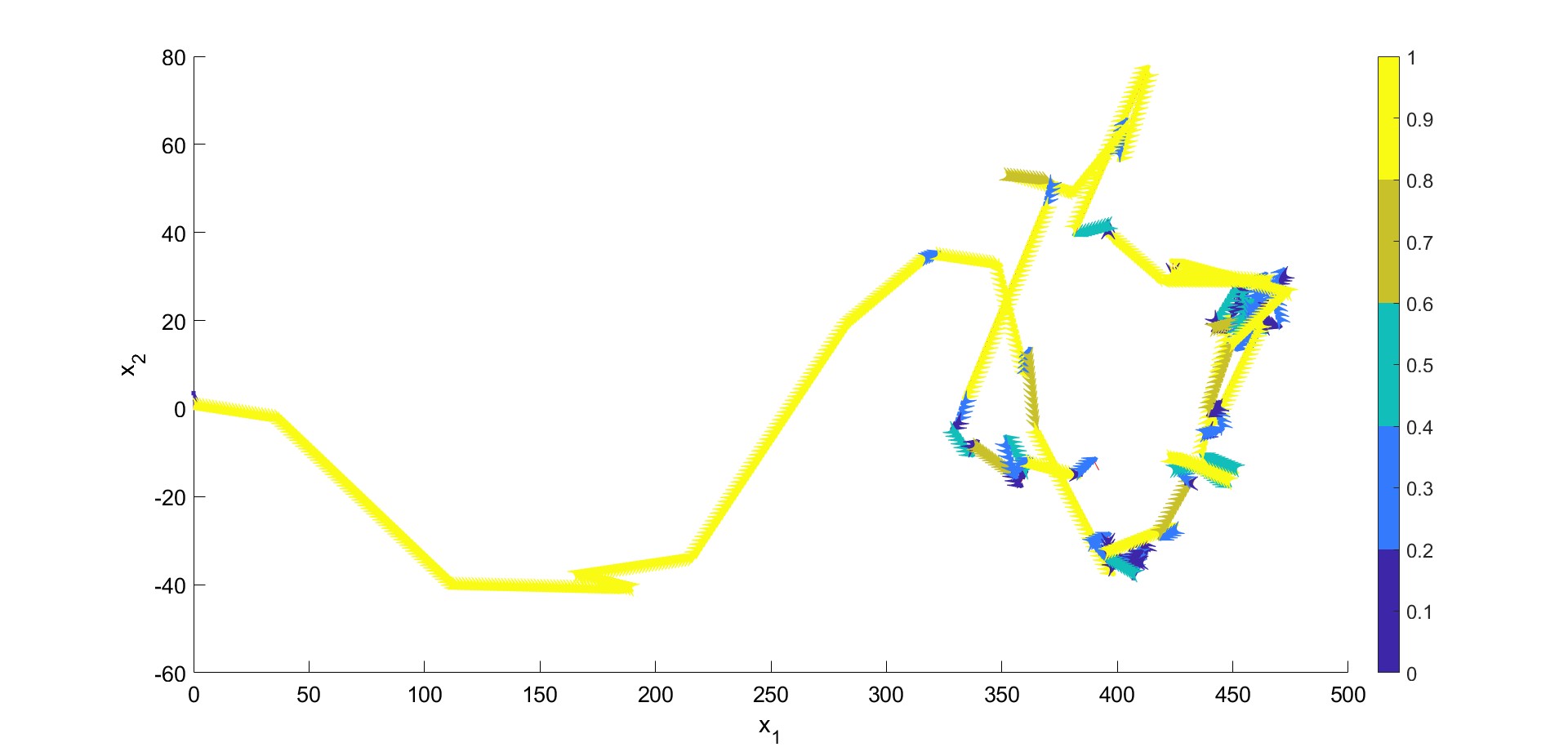}
  \includegraphics[width=13cm]{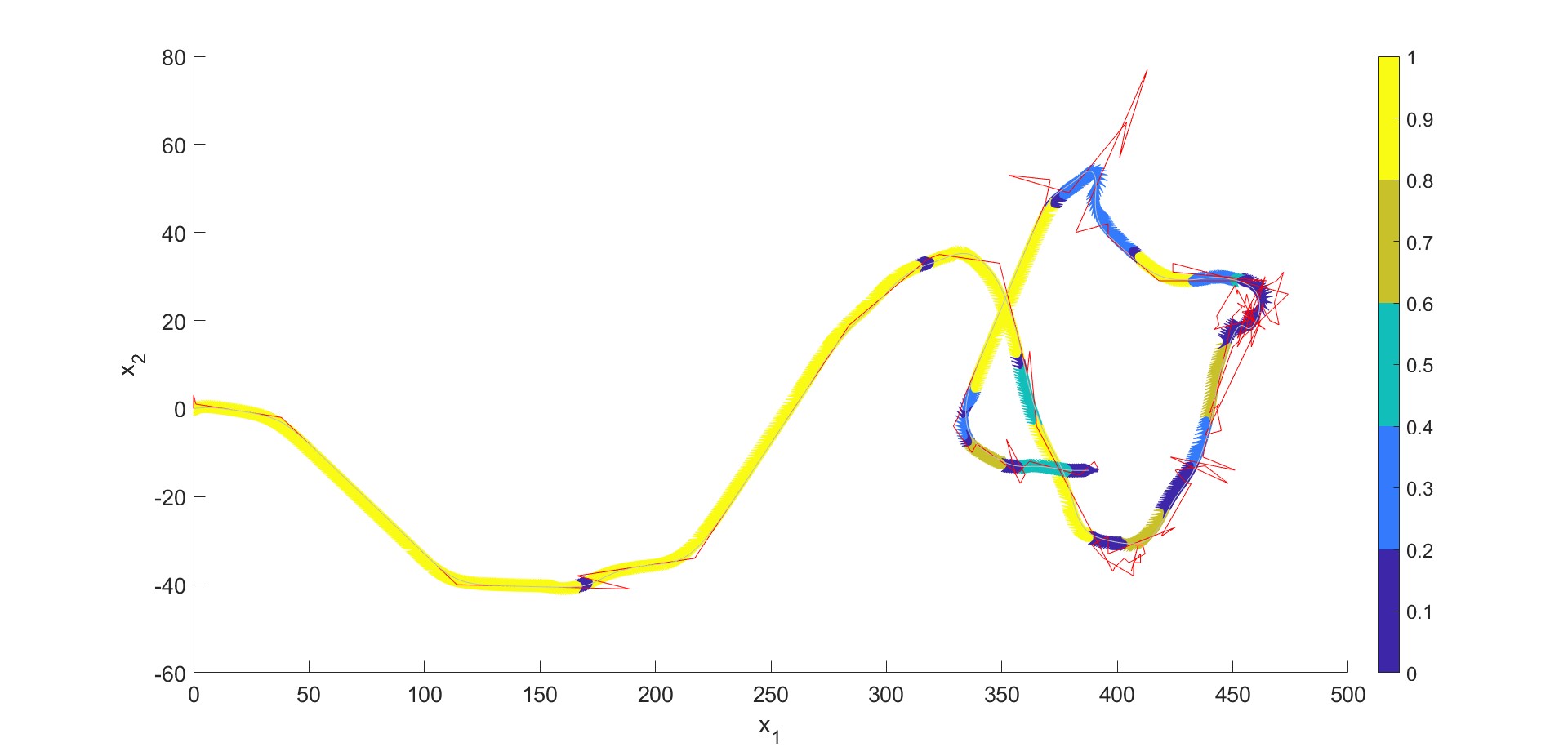}
  \qquad
  \caption{Top: initial trajectory and plotting of the velocity vectors. Bottom: initial trajectory (red line), evolved curve, and plotting of the velocity vectors. The vectors are normalized for better visualization, color represents the norm.}
  \label{Fig10}
\end{figure}

Let us now focus on the disappearing segments. We add the following condition: 
 if $L_{j}^{m+1}<h_{min}$, then $L_{j}^{m+1}=0$, where
\[
 h_{min}=\min_{1\leq i\leq n+1} h_i.
\]
If we obtain $L_{j}^{m+1}=0$ that means the first and the last endpoints of that segment at time step index $m+1$ became the same point, then $\mathcal{I}(u_{j})=\mathcal{I}(u_{j-1})$. So we modify the algorithm for finding the new indexes $i$ for $\mathcal{I}(u_j)$ adding the following condition
\begin{itemize}
 \item if $r_{j}^{m+1}=\frac{L_{j}^{m+1}}{L_s}=0$ then $\mathcal{I}(u_j)=\mathcal{I}(u_{j-1})$.
\end{itemize}
This condition will cause that, when we apply formula (\ref{EQ3.36}) for every time step index $m+k$, with $k>0$, we will obtain
\[
 L_{j}^{m+k}=L_{j}^{m+1}=0.
\]
\\
To get the new velocity estimation consider $(\Delta_t)_j$ the time between the two endpoints $\textbf{x}(u_{j-1})$ and $\textbf{x}(u_j)$ of the segment $j$. In the beginning, every segment has the same value of $(\Delta_t)_j$; later, if a segment does not disappear, we consider the same $(\Delta_t)_j$ as in the beginning. On the other hand, if a segment  $j$ disappears, we add $\frac{(\Delta_t)_j}{2}$ to the first previous not disappeared segment,
\[
 (\Delta_t)_{j-1}=(\Delta_t)_{j-1}+\frac{(\Delta_t)_j}{2},
\]
 and $\frac{(\Delta_t)_j}{2}$ to the first following not disappeared segment
 \[
 (\Delta_t)_{j+1}=(\Delta_t)_{j+1}+\frac{(\Delta_t)_j}{2}.
\]
To find the new velocity estimation for the grid points on the evolved curve, we apply the following algorithm. Set $j=1$, $i=1$ and repeat until $j=M-1$:
\begin{itemize}
 \item If $r_{j}^{m+1}> 0$ 
 \begin{equation}
  \lvert \textbf{v} \rvert=\frac{L^{d}_{j}}{(\Delta_t)_j}.
 \end{equation}
 \item for $i=\mathcal{I}(u_{j-1}),...,\mathcal{I}(u_j)$
\begin{equation}
 \textbf{v}_{i}^{m+1}=\lvert \textbf{v} \rvert\frac{(\textbf{x}_{i}^{m+1}-\textbf{x}_{i-1}^{m+1})}{\lvert \textbf{x}_{i}^{m+1}-\textbf{x}_{i-1}^{m+1}\rvert}.
\end{equation}
\end{itemize}
If the segment does not disappear, we consider the grid points inside that segment to have constant velocity equal to the new length of the segment divided by the time interval $(\Delta_t)_j$. 
On the other hand, if the $j$-th segment disappears, the time $(\Delta_t)_j$ is distributed between the previous and following segments that have not disappeared.\\
For the numerical experiments, to first check our method, we set as initial condition a semi-ellipse and we considered the evolution of the curve with $\lambda=0$, $\delta=0.05$, $\omega=1$, and $\tau=0.001$. We calculated the velocities using the algorithm described above: the result is shown in Fig. \ref{Fig7} after $1000$ time steps. Once we got the velocities for every grid point, we normalized the length of the vectors to obtain better visualization: the color represents the norm. Then, if the arrow is yellow it means in that grid point the velocity is high. On the other side, if the arrow is blue or light blue, at that point the velocity is slower. As we expected, the length of the segment in the center of the discretized semi-ellipse decreases faster than the others. Indeed, that is the part of the ellipse with the highest curvature.
Then, we focused on the real trajectories of the macrophages. We considered the evolution of the trajectories with $\lambda=1$, $\delta=0.005$, $\omega=1$, and $\tau=0.0001$. Results for three different trajectories are shown in Figs. \ref{Fig8}-\ref{Fig10}. We considered the evolved curve obtained using the mean Hausdorff distance as a stopping criterion as described in the last paragraph of section \ref{S3.1}. We observe that the velocity is slowed down in the random parts of the original trajectory that are smoothed during the evolution.  
\begin{section}{Acknowledgement}
 This work has received funding from the European Union's Horizon 2020 research and innovation programme under the Marie Sklodowska-Curie grant agreement No 955576 and by the grants APVV-19-0460, VEGA 1/0436/20.
\end{section}

\bibliographystyle{apa}

\end{document}